\documentclass[a4paper,11pt]{article}
\usepackage{adjustbox}
\usepackage{aligned-overset}
\usepackage{amsmath,amsthm}
\usepackage{authblk}
\usepackage[style=numeric-comp,maxbibnames=99,bibencoding=utf8,giveninits=true,backend=biber]{biblatex}
\usepackage{booktabs}
\usepackage{longtable}
\usepackage{multirow}
\usepackage{cancel}
\usepackage{cases}
\usepackage{colortbl}
\usepackage[hmargin={26mm,26mm},vmargin={30mm,35mm}]{geometry}
\usepackage{nicefrac}
\usepackage[colorlinks,allcolors={blue}]{hyperref}
\usepackage{paralist}
\usepackage{subcaption}
\usepackage{enumitem}
\usepackage[compact,small]{titlesec}
\usepackage{tikz-cd}
\usetikzlibrary{arrows}
\usepackage{multicol}
\usepackage{comment}

\usepackage{newtxtext}
\usepackage{newtxmath}

\bibliography{ddr_bgg}

\newcommand{\email}[1]{\href{mailto:#1}{#1}}

\numberwithin{equation}{section}

\newtheorem{theorem}{Theorem}

\newtheorem{lemma}[theorem]{Lemma}

\theoremstyle{remark}
\newtheorem{remark}[theorem]{Remark}
\theoremstyle{definition}

\newcommand{\st}{\,:\,}
\newcommand{\Real}{\mathbb{R}}

\newcommand{\Zintegers}{\mathbb{Z}}

\DeclareRobustCommand{\bvec}[1]{\boldsymbol{#1}}
\newcommand{\uvec}[1]{\underline{\bvec{#1}}}
\newcommand{\cvec}[1]{\bvec{\mathcal{#1}}}

\newcommand{\Cspace}[1]{\mathcal C_{#1}}

\DeclareMathOperator{\GRAD}{\bf grad}
\DeclareMathOperator{\CURL}{\bf curl}
\DeclareMathOperator{\DIV}{div}

\DeclareMathOperator{\ROT}{rot}
\DeclareMathOperator{\VROT}{\bf rot}
\DeclareMathOperator{\hess}{\boldsymbol{\mathcal{H}}}

\newcommand{\compl}{{\rm c}}

\newcommand{\symbolproj}{\pi}
\newcommand{\lproj}[2]{\symbolproj_{\mathcal{P},#2}^{#1}}
\newcommand{\vlproj}[2]{\boldsymbol{\symbolproj}_{\cvec{P},#2}^{#1}}
\newcommand{\Rproj}[2]{\bvec{\symbolproj}_{\cvec{R},#2}^{#1}}

\newcommand{\Xgrad}[2]{\underline{X}_{\GRAD,#2}^{#1}}
\newcommand{\Xcurl}[2]{\underline{\bvec{X}}_{\CURL,#2}^{#1}}

\newcommand{\Igrad}[2]{\underline{I}_{\GRAD,#2}^{#1}}

\newcommand{\GE}[1]{G_E^{#1}}

\newcommand{\uGT}[1]{\uvec{G}_T^{#1}}
\newcommand{\uGh}[1]{\uvec{G}_h^{#1}}

\newcommand{\uCh}[1]{C_h^{#1}}

\newcommand{\edges}[1]{\mathcal{E}_{#1}}
\newcommand{\vertices}[1]{\mathcal{V}_{#1}}

\newcommand{\ET}{\edges{T}}

\newcommand{\VT}{\vertices{T}}

\newcommand{\VE}{\vertices{E}}

\newcommand{\normal}{\bvec{n}}
\newcommand{\tangent}{\bvec{t}}

\newcommand{\letterPoly}{\mathcal{P}}
\newcommand{\Poly}[1]{\letterPoly^{#1}}
\newcommand{\Roly}[1]{\cvec{R}^{#1}}
\newcommand{\cRoly}[1]{\cvec{R}^{\compl,#1}}

\newcommand{\stokes}{2}
\newcommand{\ddr}{1}

\newcommand{\vvvert}{\vert\kern-0.25ex\vert\kern-0.25ex\vert}

\DeclareMathOperator{\Ker}{Ker}
\DeclareMathOperator{\Image}{Im}

\newcommand{\Th}{\mathcal{T}_h}
\newcommand{\Eh}{\mathcal{E}_h}
\newcommand{\Vh}{\mathcal{V}_h}

\newcommand{\Prot}[2]{\bvec{P}_{\VROT,#2}^{#1}}

\newcommand{\XSgrad}[1]{\underline{H}_{2}^{k}(#1)}

\newcommand{\XSrot}[1]{\underline{H}_{1}^{k+1}(#1)^2}

\newcommand{\ISgrad}{\underline{I}_{\stokes,h}^{k}}
\newcommand{\ISgradT}{\underline{I}_{\stokes,T}^{k}}

\newcommand{\SGRAD}[1]{\underline{\bvec{G}}^{k}_{2,#1}}
\newcommand{\SROT}[1]{R^k_{1,#1}}

\newcommand{\Gqen}{G^{\normal}_{q,E}}
\newcommand{\Gqv}{\bvec G_{q,V}}
\newcommand{\Gqet}{G^{\tangent}_{2,E}}
\newcommand{\nablaT}{\bvec{G}_{2,T}^{k-1}}

\newcommand{\Id}{\mathrm{Id}}
\newcommand{\sskw}{\mathop{\mathrm{sskw}}}
\newcommand{\sskwh}{\sskw\nolimits_h}

\newcommand{\XdRgrad}[1]{\underline{H}_{1}^{k+1}(#1)}

\newcommand{\XdRrot}[1]{\underline{\boldsymbol{H}}_{\operatorname{\boldsymbol{\mathrm{rot}}}}^{k+1}(#1)}

\newcommand{\tXdRgrad}[1]{\XdRgrad{#1}^2}

\newcommand{\tXdRrot}[1]{\XdRrot{#1}}

\newcommand{\tGRAD}[1]{\underline{\bvec G}_{1,#1}^{k+1}}
\newcommand{\tROT}[1]{\bvec{R}_{\operatorname{\boldsymbol{\mathrm{rot}}},#1}^{k+1}}
\newcommand{\IdRgrad}{\underline{\bvec{I}}_{\ddr,h}^{k+1}}
\newcommand{\IdRgradT}{\underline{\bvec{I}}_{\ddr,T}^{k+1}}
\newcommand{\IdRrot}{\underline{\bvec{I}}_{\operatorname{\boldsymbol{\mathrm{rot}}},h}^{k+1}}
\newcommand{\IdRrotT}{\underline{\bvec{I}}_{\operatorname{\boldsymbol{\mathrm{rot}}},T}^{k+1}}

\newcommand{\tGT}{\bvec{G}_{1,T}^{k}}
\newcommand{\tGE}{\bvec{G}^{k+1}_{1,E}}

\newcommand{\bs}{{\scriptscriptstyle \bullet}}

\newcommand{\XdRrotS}[1]{\underline{\boldsymbol{H}}_{\operatorname{\boldsymbol{\mathrm{rot}}}}^{k+1}(#1,\mathbb{S})}

\newcommand{\Hess}[1]{\underline{\bvec{\mathcal{H}}}^{k+1}_{#1}}

\newcommand{\redex}[1]{\mathfrak{#1}}

\newcommand{\Egradh}{\underline{\redex{E}}_{\GRAD,h}}
\newcommand{\EgradT}{\underline{\redex{E}}_{\GRAD,T}}
\newcommand{\EgradE}{\underline{\redex{E}}_{\GRAD,E}}
\newcommand{\Egrad}[2]{\redex{E}_{\letterPoly,#2}^{#1}}
\newcommand{\Rgradh}{\underline{\redex{R}}_{\GRAD,h}}

\newcommand{\Eroth}{\underline{\bvec{\redex{E}}}_{\ROT,h}}
\newcommand{\ErotT}{\underline{\bvec{\redex{E}}}_{\ROT,T}}

\newcommand{\Erot}[2]{\boldsymbol{\redex{E}}_{\boldsymbol{R},#2}^{#1}}
\newcommand{\Rroth}{\underline{\redex{R}}_{\ROT,h}}
\newcommand{\RrotT}{\underline{\redex{R}}_{\ROT,T}}
\newcommand{\dof}[2]{\underline{#1}_{#2}}
\newcommand{\dofh}[1]{\underline{#1}_h}
\newcommand{\dofT}[1]{\underline{#1}_T}

\usepackage{marginnote}

\usepackage{pgfplots,pgfplotstable}
\usepackage{tikz-cd}

\graphicspath{{figures/}}

\newcommand{\drawPolygon}[7]{%
  \def\nsides{#1}%
  \def\radius{#2}%
  \def\showNodal{#3}%
  \def\showGradients{#4}%
  \def\pEdge{#5}%
  \def\varsEdge{1}%
  \def\pNormal{#6}%
  \def\varsNormal{1}%
  \def\nInterior{#7}%

  \foreach \i in {1,...,\nsides} {
    \pgfmathsetmacro{\angle}{360/\nsides*(\i - 1)}
    \pgfmathsetmacro{\x}{\radius*cos(\angle)}
    \pgfmathsetmacro{\y}{\radius*sin(\angle)}
    \coordinate (P\i) at (\x,\y);
  }

  \draw[thick] (P1)
  \foreach \i in {2,...,\nsides} { -- (P\i) }
  -- cycle;

  \ifnum\showNodal=1
  \foreach \i in {1,...,\nsides} {
    \fill (P\i) circle (2pt);
  }
  \fi

  \ifnum\showGradients=1
  \foreach \i in {1,...,\nsides} {
    \draw (P\i) circle (4pt);
  }
  \fi

  \ifnum\pEdge>-1
  \pgfmathtruncatemacro{\nEdge}{int((\pEdge+\varsEdge)!/(\pEdge!*\varsEdge!))}
  \foreach \i in {1,...,\nsides}{
    \pgfmathtruncatemacro{\j}{mod(\i,\nsides)+1}%
    \foreach \k in {1,...,\nEdge}{
      \pgfmathsetmacro{\t}{\k/(\nEdge+1)}
      \path (P\i) -- (P\j) coordinate[pos=\t] (E\i\k);
      \fill (E\i\k) circle (1.5pt);
    }
  }
  \fi

  \ifnum\pNormal>-1
  \pgfmathtruncatemacro{\nNorm}{int((\pNormal+\varsNormal)!/(\pNormal!*\varsNormal!))}
  \foreach \i in {1,...,\nsides}{
    \pgfmathtruncatemacro{\j}{mod(\i,\nsides)+1}
    \pgfmathsetmacro{\angleA}{360/\nsides*(\i - 1)}
    \pgfmathsetmacro{\angleB}{360/\nsides*(\j - 1)}
    \pgfmathsetmacro{\edgeang}{atan2(sin(\angleB)-sin(\angleA),cos(\angleB)-cos(\angleA))}
    \pgfmathsetmacro{\normang}{\edgeang-90}
    \foreach \k in {1,...,\nNorm}{
      \pgfmathsetmacro{\t}{\k/(\nNorm+1)}
      \path (P\i) -- (P\j) coordinate[pos=\t] (N\i\k);
      \draw[-{Latex[length=4pt]}, thick]
      (N\i\k) -- ++({0.6*cos(\normang)},{0.6*sin(\normang)});
    }
  }
  \fi

  \ifnum\nInterior>0
  \ifnum\nInterior=1
  \fill (0,0) circle (1.5pt);
  \else
  \pgfmathsetmacro{\rInt}{0.35*\radius}
  \foreach \k in {1,...,\nInterior} {
    \pgfmathsetmacro{\theta}{360/\nInterior*(\k - 1)}
    \pgfmathsetmacro{\xi}{\rInt*cos(\theta)}
    \pgfmathsetmacro{\yi}{\rInt*sin(\theta)}
    \fill (\xi,\yi) circle (1.5pt);
  }
  \fi
  \fi
}

\newcommand{\DrawHtwo}{%
  \tikz[baseline=(P1.base), every node/.style={inner sep=0}]{
    \drawPolygon{6}{0.6}{1}{1}{-1}{0}{-1}%
  }%
}

\newcommand{\DrawHone}{%
  \tikz[baseline=(P1.base), every node/.style={inner sep=0}]{
    \drawPolygon{6}{0.6}{1}{0}{0}{-1}{-1}%
  }%
}

\newcommand{\DrawHrot}{%
  \tikz[baseline=(P1.base), every node/.style={inner sep=0}]{
    \drawPolygon{6}{0.6}{0}{0}{1}{-1}{2}%
  }%
}

\newcommand{\DrawPzero}{%
  \tikz[baseline=(P1.base), every node/.style={inner sep=0}]{
    \drawPolygon{6}{0.6}{0}{0}{-1}{-1}{1}%
  }%
}

\newcommand{\DrawPone}{%
  \tikz[baseline=(P1.base), every node/.style={inner sep=0}]{
    \drawPolygon{6}{0.6}{0}{0}{-1}{-1}{3}%
  }%
}

\begin{document}

\title{Design and homological analysis of twisted and BGG Stokes-de Rham complexes on polygonal meshes}

\author[1]{Daniele A. Di Pietro}
\author[1,2]{J\'er\^ome Droniou}
\author[3]{Kaibo Hu}
\author[1]{Arax Leroy}
\affil[1]{IMAG, Univ. Montpellier, CNRS, Montpellier, France, \email{daniele.di-pietro@umontpellier.fr}, \email{jerome.droniou@cnrs.fr}, \email{arax.leroy@umontpellier.fr}}
\affil[2]{School of Mathematics, Monash University, Melbourne, Australia}
\affil[3]{Mathematical Institute, University of Oxford, UK. \email{kaibo.hu@maths.ox.ac.uk}}
\maketitle

\begin{abstract}
  We design a discrete Bernstein--Gelfand--Gelfand (BGG) diagram in a two-dimensional setting on polygonal meshes based on the DDR framework; the diagram is made of a  discrete Stokes polygonal complex and a tensorised Discrete de Rham complex, and the BGG construction leads to novel twisted and Hessian complexes applicable on generic polygonal meshes; such complexes enable, in particular, discretisations of plate problems. Complete homological properties of the discrete Stokes complex and of the complexes built from the BGG diagram are established.
\end{abstract}


\section{Introduction}\label{sec:Introduction}

Hilbert complexes naturally arise in the analysis of important classes of partial differential equations encountered across various areas of physics and engineering \cite{Arnold.Falk.ea:10}. In this work, we focus on the two-dimensional Stokes and Hessian complexes, which are respectively relevant to incompressible fluid mechanics and plate theory. It is known \cite{Arnold.Hu:21,Cap.Hu:24} that these complexes, along with the standard de Rham complex, are linked through the Bernstein--Gelfand–Gelfand (BGG) construction \cite{Cap.Slovak.ea:01,Bernstein.Gelfand.ea:75}, a framework that has already been exploited to derive finite element complexes (see, e.g., \cite{Christiansen.Hu.ea:18,Chen.Huang:24,Bonizzoni.Hu.ea:25}). Here, we investigate for the first time its use in constructing discrete Stokes and Hessian complexes on general polygonal meshes.
\smallskip

Denote by $\Omega$ a two-dimensional polygonal domain.
We consider the following BGG diagram, which stacks the \emph{Stokes complex} (i.e., a version of the de Rham complex with increased regularity) above a tensorised version of the usual de Rham complex:
\begin{equation}\label{eq:double.complex.cont}
  \begin{tikzcd}[column sep=2.5em]
    \text{Stokes:}
    & 0 \arrow{r}{}
    & H_2(\Omega)\arrow{r}{\GRAD}
    & H_1(\Omega)^2\arrow{r}{\ROT}
    & L^2(\Omega)\arrow{r}{}
    & 0 \\
    \text{de Rham:}
    & 0 \arrow{r}{}
    & H_1(\Omega)^2\arrow{r}{\GRAD}\arrow[leftrightarrow]{ur}{\Id}
    & \boldsymbol{H}_{\VROT}(\Omega)\arrow{r}{\VROT}\arrow{ur}{\sskw}
    & L^2(\Omega)^2\arrow{r}{}
    & 0.
  \end{tikzcd}
\end{equation}
Above, $L^2(\Omega)$ denotes the space of square-integrable functions on $\Omega$ while,
for any integer $m \ge 0$, $H_m(\Omega)$ is the usual Hilbert space spanned by functions that have partial derivatives up to degree $m$ in $L^2(\Omega)$.
Additionally, for smooth enough scalar-, vector-, and tensor-valued functions
\[
q : \Omega \to \Real\,,\quad \bvec{v} = \begin{pmatrix}
  v_1 \\ v_2
\end{pmatrix}: \Omega \to \Real^2\quad\text{ and }\quad
\bvec{\tau} = \begin{pmatrix}
  \tau_{11} & \tau_{12} \\ \tau_{21} & \tau_{22}
\end{pmatrix}: \Omega \to \Real^{2\times 2},
\]
we have set
\[
\GRAD q \coloneq \begin{pmatrix} \partial_1 q \\ \partial_2 q
\end{pmatrix},\quad
\ROT \bvec{v} \coloneq \partial_1 v_2 - \partial_2 v_1,\quad
\GRAD \bvec{v} \coloneq \begin{pmatrix}
  \partial_1 {v}_1 & \partial_2 {v}_1 \\
  \partial_1 {v}_2 & \partial_2 {v}_2
\end{pmatrix},\quad
\VROT \bvec{\tau} \coloneq \begin{pmatrix}
  \partial_1 \tau_{12} - \partial_2 \tau_{11} \\
  \partial_1 \tau_{22} - \partial_2 \tau_{21} \\
\end{pmatrix}
\]
and
\begin{equation}\label{eq:sskw}
  \sskw \bvec{\tau} \coloneq \tau_{12} - \tau_{21}.
\end{equation}
Finally, $\boldsymbol{H}_{\VROT}(\Omega) \coloneq \left\{ \bvec{v} \in L^2(\Omega)^2 \st \ROT \bvec{v} \in L^2(\Omega) \right\}$.
We will show in Section \ref{sec:bgg} that, starting from the BGG diagram \eqref{eq:double.complex.cont}, one can derive the following \emph{Hessian complex}:
\begin{equation}\label{BGG-Hessian-HD}
  \begin{tikzcd}
    0 \arrow{r} &H_{2}(\Omega) \arrow{r}{\hess} & \boldsymbol H_{\boldsymbol \ROT}(\Omega, \mathbb{S})\arrow{r}{\boldsymbol \ROT} &L^{2}(\Omega)^{2}  \arrow{r}{} &0,
  \end{tikzcd}
\end{equation}
where $\hess \coloneqq\GRAD \GRAD$ is the Hessian operator
and, denoting by $\mathbb{S}$ the space of $2\times 2$ symmetric matrices and by $L^2(\Omega,\mathbb{S})$ the space of square-integrable functions $\Omega\to\mathbb{S}$, we have set
\[
\boldsymbol H_{\VROT}(\Omega, \mathbb{S}) \coloneq \left\{ \boldsymbol{\sigma} \in L^2(\Omega,\mathbb{S}): \VROT \boldsymbol{\sigma} \in L^2(\Omega)^2 \right\}.
\]
The Hessian complex is relevant, e.g., in the analysis of the Kirchhoff--Love plates equation;
see Section \ref{sec:complex.in.elasticity} for an overview.

\smallskip

The design of discrete Hilbert complexes has been a central topic in finite element research over the past two decades, with several fundamental questions still unresolved.
In what follows, we will briefly review a few contributions relevant to the present work.
Around 2000, the construction of stable finite element pairs for elasticity equations in mixed form using polynomial shape functions saw a breakthrough with the introduction of the Arnold--Winther element \cite{Arnold.Winther:02}.
This element is part of a discrete elasticity complex that begins with the Argyris ($\Cspace{1}$) space.
This complex can be incorporated into a BGG diagram alongside finite element de Rham complexes \cite{Arnold.Falk.ea:06}, providing a BGG interpretation of the construction.

Another notable finite element elasticity pair is the Hu--Zhang element \cite{Hu.Zhang:15}.
Similar to the Arnold--Winther element, this element is also part of a complex that starts with the Argyris space.
Furthermore, the Hu--Zhang element and its associated complex can be derived using a discrete version of the BGG diagram \cite{Christiansen.Hu.ea:18} (see Section~\ref{sec:fe.comparison} below for further details).

Subsequently, several other examples have emerged where finite element Stokes and de Rham complexes are integrated into BGG diagrams to develop finite elements for the Hellinger--Reissner principle and elasticity complexes. In particular, various scalar $\Cspace{1}$ spline spaces can serve as starting points \cite{Lai.Schumaker:07}. For instance, by beginning with the Hsieh--Clough--Tocher macroelement on the Alfeld (Clough--Tocher) split, where a triangle is divided into three subtriangles, one can derive a finite element elasticity complex on the same triangulation \cite{Christiansen.Hu:23}.
This elasticity pair can be seen as a generalization of the earlier construction by Johnson and Mercier \cite{Johnson.Mercier:78}. In the corresponding BGG construction, the first row of the diagram represents a Stokes complex as described in \cite{Christiansen.Hu:18}, with the Stokes pair originally constructed by Arnold and Qin \cite{Arnold.Qin:92}. The second row consists of a standard finite element de Rham complex defined on the Alfeld split.
A similar construction for criss-cross grids was presented in \cite{Hu:22}.

Finally, constructions of complexes containing second-order derivatives have also been carried out without resorting to BGG diagrams; see, e.g., \cite{Chen.Huang:22,Chen.Huang:22*1,Chen.Huang:22*2,Hu.Liang.ea:22,Hu.Lin.ea:23,Hu.Lin.ea:25}.
These constructions, however, do not expose the entire physics contained in BGG diagrams such as, e.g., twisted complexes.
\smallskip

In the present work, we apply for the first time the BGG construction to polytopal approximations of Hilbert complexes by designing a discrete counterpart of the BGG diagram \eqref{eq:double.complex.cont}. Unlike standard finite elements, polytopal methods are formulated on general meshes, allowing for elements of arbitrary shape and hanging nodes \cite{Di-Pietro.Droniou:20}. The derivation of polytopal Hilbert complexes is a relatively recent line of research; see \cite{Beirao-da-Veiga.Mora.ea:19,Di-Pietro.Droniou.ea:20,Di-Pietro.Droniou:23,Di-Pietro.Droniou:23*1,Di-Pietro.Droniou.ea:23,Hanot:23,Di-Pietro:24} and references therein.

Within this framework, the discrete counterpart of the bottom sequence in \eqref{eq:double.complex.cont} is chosen as the serendipity Discrete de Rham (DDR) complex introduced in \cite{Di-Pietro.Droniou:23*1}. For the top sequence, we construct a new Stokes complex by reverse-engineering the head and tail spaces from the middle one. This construction differs from the discrete Stokes (DS) complex proposed in \cite{Hanot:23}. In particular, our choice of the bottom sequence yields discrete spaces whose components all belong to full polynomial spaces, resulting in a simpler and more cost-effective implementation (see Remark~\ref{rem:serendipity}).

Building on existing results for the DDR method \cite{Di-Pietro.Droniou.ea:23,Di-Pietro.Droniou:23*1,Bonaldi.Di-Pietro.ea:25}, we prove that the cohomology of the resulting DS complex is isomorphic to that of the continuous de Rham complex. From the discrete BGG diagram, we then derive a new polygonal discrete Hessian complex by applying the BGG procedure. We further show that this complex has the same cohomology as its continuous counterpart, even on domains with nontrivial topology. While homological properties of discrete BGG complexes are often established only in the case of trivial topology (i.e., exactness), we note that our dimension-counting argument can be adapted to a broader class of discrete complexes.

\smallskip

The rest of this work is organised as follows.
In Section~\ref{sec:bgg} we briefly recall the BGG construction and apply it to the derivation of the twisted and BGG complexes associated to Diagram~\eqref{eq:double.complex.cont}.
After describing the discrete setting in Section~\ref{sec:setting}, in Section~\ref{sec:discrete.bgg} we recall the construction of the DDR complex and describe the new Discrete Stokes complex, and show how they are organised in a BGG diagram.
The homological properties of the DS complex are studied in Section~\ref{sec:homological_properties}.
In Section \ref{sec:twisted_and_bgg_complexes}, we apply the BGG machinery to obtain the corresponding twisted and BGG polygonal complexes, the latter being a discrete version of the two-dimensional Hessian complex. For ease of reference, we gathered in an appendix some notations of spaces and operators that are used throughout the article.


\section{Twisted and BGG complexes}\label{sec:bgg}

\subsection{Basic principles of the BGG construction}

We briefly recall the BGG construction of \cite{Arnold.Hu:21}.
A \emph{BGG diagram}
\begin{equation}\label{BGG-diagram}
  \begin{tikzcd}[column sep=2em, row sep=3em]
    {\cdots} \arrow{r} &V^{k-2}  \arrow{r}{d^{k-2}} &[3em]V^{k-1} \arrow{r}{d^{k-1}} &[3em]V^{k} \arrow{r}{d^{k}} &[3em] V^{k+1} \arrow{r}{} &\cdots\\
    \cdots \arrow{r}& W^{k-2}\arrow{r}{d^{k-2}} \arrow[ur, "{S^{k-2}}"]&{ W^{k-1} } \arrow{r}{{d^{k-1}}} \arrow[ur, "{S^{k-1}}"]&{  W^{k}} \arrow{r}{{d^{k}}}\arrow[ur, "{S^{k}}"] &{W^{k+1} }\arrow{r}{} &\cdots
  \end{tikzcd}
\end{equation}
consists of complexes connected by algebraic operators $S^{\bs}$ in an anti-commuting diagram satisfying
\begin{equation}\label{anti-commuting}
  dS=-Sd.
\end{equation}
Here, $V^{i}$ and $W^{i}$ are Hilbert spaces and $d^{i}$ are linear operators. Typically, each row is a scalar- or vector-valued de~Rham complex.

From \eqref{BGG-diagram}, we can immediately read out the
\emph{twisted complex}:
\begin{equation}\label{cplx:twisted}
  \begin{tikzcd}[ampersand replacement=\&, column sep=2em]
    \cdots\arrow{r}\&
    \begin{pmatrix}
      V^{k-1} \\
      W^{k-1}
    \end{pmatrix}
    \arrow{r}{
      \begin{pmatrix} d^{k-1} & -S^{k-1} \\ 0 & d^{k-1} \end{pmatrix}
    }\&[4em]  \begin{pmatrix}
      V^{k}  \\
      W^{k}
    \end{pmatrix} \arrow{r}{
      \begin{pmatrix} d^{k} & -S^{k} \\ 0 & d^{k} \end{pmatrix}
    } \&[4em] \begin{pmatrix}
      V^{k+1}   \\
      W^{k+1}
    \end{pmatrix} \arrow{r}{} \&\cdots,
  \end{tikzcd}
\end{equation}
with operators $A^{i}\coloneqq  \begin{pmatrix} d^{i} & -S^{i} \\ 0 & d^{i} \end{pmatrix}$.
The sequence \eqref{cplx:twisted} is a complex, i.e., $A^{k+1}\circ A^{k}=0$ for any $k$, thanks to the anti-commutativity \eqref{anti-commuting}.   By \cite[Lemma 6]{Arnold.Hu:21}, the dimension of the cohomology of \eqref{cplx:twisted} is less than or equal to the sum of the cohomology dimensions of the inputs $(V^{\bs}, d^{\bs})$ and $(W^{\bs}, d^{\bs})$. Equality holds, i.e., the cohomology of the output is isomorphic to the input, if and only if $S^{\bs}$ induces zero maps on cohomology \cite[Lemma 7]{Arnold.Hu:21}. A typical assumption leading to this property is the existence, for all $i$, of $K^{i}:W^{i}\to V^{i}$ satisfying $S^{k}=d^{k}K^{k}-K^{k+1}d^{k}$. This holds for the example of the Hessian complex considered below, where $K^i(\bvec{w}):=\bvec{w}\cdot \bvec{x}$ (row-wise multiplication by the position vector) for any $\bvec{w}\in W^k$.

In typical applications, there exists an index $J$ such that
$S^k$ is injective for $k \le J$ and surjective for $k \ge J$,
with a bijective $S^{J}$ in the middle.
Then, the spaces in the diagram can be decomposed as the kernel and cokernel of the $S^{\bs}$ maps, i.e.,
\begin{equation*}
  \begin{tikzcd}[column sep=2em, row sep=1cm]
    {\cdots} \hspace*{-2em}
    & \Image(S^{J-2}) \oplus \Image(S^{J-2})^{\perp} \arrow{r}{d^{J-1}}
    & \Image(S^{J-1}) \oplus \Image(S^{J-1})^{\perp} \arrow{r}{d^{J}}
    & V^{J+1}=\Image(S^{J})
    &\hspace*{-2em} \cdots\\
    \cdots \hspace*{-2em}
    & W^{J-1} \arrow{r}{{d^{J-1}}} \arrow[ur, "{S^{J-1}}"]
    & W^{J} \arrow{r}{{d^{J}}}\arrow[ur, "{S^{J}}"]
    & \Ker(S^{J+1})\oplus \Ker(S^{J+1})^{\perp}
    &\hspace*{-2em}\cdots
  \end{tikzcd}
\end{equation*}
The following BGG complex  is obtained by eliminating components connected by $S^{\bs}$:
\begin{equation} \label{BGG-seq}
  \begin{tikzcd} [column sep=2em, row sep=1cm]
    \cdots \Image(S^{J-2})^{\perp} \ar[r,"P_{ \Image(S^{J-1})^{\perp} }\circ d^{J-1}"]  &[3.5em]    \Image(S^{J-1})^{\perp}\ar[r,"d^J"] & \rule{0em}{1em}\rule{1em}{0em} \ar{ld}[swap]{(S^{J})^{-1}} \\
    &   \rule{0em}{1em}\rule{1em}{0em}\ar[r,"d^J"] &   \Ker(S^{J+1})\ar[r,"d^{J+1}"] &  \Ker(S^{J+2})   \cdots,
  \end{tikzcd}
\end{equation}
where $P_{ \Image(S^{J-1})^{\perp} }$ denotes the  projection  onto  $\Image(S^{J-1})^\perp$.
The anti-commutativity \eqref{anti-commuting} also implies that $d^{k}$ maps $ \Ker(S^{k})$ to $ \Ker(S^{k+1})$ in \eqref{BGG-seq}.

\subsection{BGG derivation of the hessian complex}

Next, we present the example relevant to this paper, i.e., the diagram \eqref{eq:double.complex.cont}, leading to the Hessian complex. The first row of \eqref{eq:double.complex.cont} is a Stokes complex, and the second row is a de~Rham complex with $L^{2}$-based Sobolev spaces. The operator $\sskw:  \boldsymbol    H_{\boldsymbol \ROT}(\Omega)^{2}  \to L^{2}(\Omega) $ is onto since the diagram \eqref{eq:double.complex.cont} anti-commutes.
The twisted complex derived from \eqref{eq:double.complex.cont} is:
\begin{equation}\label{twisted-stokes-HD}
  \begin{tikzcd}[ampersand replacement=\&, column sep=2em]
    0\arrow{r}\&
    \begin{pmatrix}
      H_{2}(\Omega)  \\
      H_{1}(\Omega)^{2}
    \end{pmatrix}
    \arrow{r}{
      \begin{pmatrix} \GRAD & -I \\ 0 & \GRAD \end{pmatrix}
    }\&[3.5em] \begin{pmatrix}
      H_{1}(\Omega)^{2}  \\
      \boldsymbol    H_{\boldsymbol \ROT}(\Omega)^{2}
    \end{pmatrix} \arrow{r}{
      \begin{pmatrix} \ROT & -\sskw \\ 0 &\boldsymbol  \ROT \end{pmatrix}
    } \&[3.5em] \begin{pmatrix}
      L^{2}(\Omega)   \\
      L^{2}(\Omega)^{2}
    \end{pmatrix} \arrow{r}{} \&0
  \end{tikzcd}
\end{equation}
and, since $\Ker(\sskw)=\mathbb{S}$, the corresponding BGG (Hessian) complex derived from \eqref{eq:double.complex.cont} is precisely \eqref{BGG-Hessian-HD}.
A major conclusion of the BGG construction is that the cohomology of the output (twisted, BGG) complexes is isomorphic to that of the input (a sum of de~Rham complexes). In our case, we have the following.

\begin{theorem}[Cohomologies of the twisted and BGG complexes]
  The cohomologies of the twisted complex \eqref{twisted-stokes-HD} and the BGG (hessian) complex \eqref{BGG-Hessian-HD} are isomorphic to a sum of the de~Rham cohomologies $\mathcal{H}^{\bs}_{\mathrm{dR}}\otimes (\mathbb{R}\times \mathbb{R}^{2})$, where $\mathcal{H}^{\bs}_{\mathrm{dR}}$ denotes the de Rham cohomology.
\end{theorem}

\subsection{Relevance to models in linear elasticity}\label{sec:complex.in.elasticity}

One of the motivations to discretise the whole BGG diagram such as \eqref{eq:double.complex.cont}, rather than directly discretising the BGG complex such as \eqref{BGG-Hessian-HD}, is that the twisted complex encodes richer physics \cite{Cap.Hu:24}. In our example,  the Hodge Laplacian problem of the twisted complex leads to the energy functional of the Reissner--Mindlin plate
\[
\mu_{c}\|\GRAD u-\bvec{w}\|_{A}^{2}+\|\GRAD \bvec{w}\|_{C}^{2},
\]
for $u\in H_{1}(\Omega)$ and $\bvec{w}\in H_{1}(\Omega)^{2} $ with proper weighted norms $\|{\cdot}\|_{A}$ and $\|{\cdot}\|_{C}$. Here $\mu_{c}$ is a physical parameter related to the thickness of the plate, describing the coupling between the vertical displacement (bending) $u$ and the rotation (shear) $\bvec{w}$.  The limit $\mu_{c}\rightarrow \infty$ forces $\bvec{w}=\GRAD u$, and thus the energy function becomes
\begin{equation}\label{energy:KL}
  \|\hess u\|_{C}^{2}.
\end{equation}
This describes the Kirchhoff--Love plate. The energy functional \eqref{energy:KL} also corresponds to the first Hodge Laplacian problems of \eqref{BGG-Hessian-HD}. Therefore, the BGG construction can be interpreted as a cohomology-preserving elimination of the rotational degrees of freedom $\bvec{w}$ from the Reissner--Mindlin model to get the Kirchhoff--Love plate.

In two dimensions, one may replace $\GRAD$-$\ROT$ in the complexes by $\CURL$-$\DIV$. Then, the Hodge Laplacian problems of the last part of the twisted complex \eqref{twisted-stokes-HD} and the BGG complex \eqref{BGG-Hessian-HD} corresponds to the mixed form of the linear Cosserat  model and linear elasticity (Hellinger--Reissner principle), respectively \cite{Cap.Hu:24,Boon.Duran.ea:25}.

\begin{remark}[Hodge Laplacian]
  Precisely defining Hodge Laplacian problems requires the language of Hilbert complexes \cite{Arnold.Falk.ea:10,Arnold:18}, which is feasible for the examples above. In the previous discussion, we omitted the details to avoid technicalities which are not straightforwardly related to the topic of this paper, i.e., the construction of discrete diagrams and complexes. In this simplified presentation, the Hodge Laplacian can be understood in a formal way as $dd^{\ast}+d^{\ast}d$, where the $d^{\ast}$ is the {\it formal} adjoint.
\end{remark}


\section{Discrete setting}\label{sec:setting}

\subsection{Mesh}

Given a two-dimensional polygonal domain $\Omega \subset \mathbb{R}^2$, we consider a polygonal mesh $\mathcal{M}_h = (\Th, \Eh, \Vh)$ of $\Omega$, where
$\Th$ is a finite collection of open, disjoint polygonal elements $T$ with diameter $h_T$, such that $\overline{\Omega} = \bigcup_{T \in \Th} \overline{T}$ and $h = \max_{T \in \Th} h_T > 0$;
$\Eh$ is a finite collection of open straight edges $E$ of length $h_E$;
$\Vh$ is the set of vertices $V$ with coordinate vector $\bvec{x}_V$, corresponding to the endpoints of the edges in $\Eh$.
Furthermore, we assume that the pair $(\Th,\Eh)$ satisfies \cite[Definition~1.4]{Di-Pietro.Droniou:20}.
In particular, each edge is contained in the boundary of at least one element, and the boundary of each element is the union of (closures of) the edges collected in the set $\ET$.
This broad definition permits, for instance, a flat portion of an element's boundary to be subdivided into several mesh edges, a situation encountered in non-conforming local mesh refinement.

For all $Y \in \Th \cup \Eh$, $\mathcal{V}_Y$ denotes the set of vertices of $Y$.
Each edge $E\in\Eh$ is endowed with an orientation determined by a fixed unit tangent vector $\tangent_E$; we then choose the unit normal $\normal_E$ such that $(\tangent_E,\normal_E)$ forms a right-handed system of coordinates.
For $T\in\Th$ and $E\in\ET$, we let $\omega_{TE} \in \{ -1, 1 \}$ be such that $\omega_{TE}\normal_E$ is the unit vector normal to $E$ pointing out of $T$.
For $E\in\Eh$ and $V\in\mathcal{V}_E$, we set $\omega_{EV}=1$ if $\tangent_E$ points in the direction of $V$, and $\omega_{EV}=-1$ otherwise.

We note that, for all $T \in \Th$, all $V \in \mathcal{V}_T$, and any family of vertex values $(\varphi_V)_{V \in \mathcal{V}_T}\in\Real^{\VT}$, we have
\begin{equation}\label{eq:sum.omegaTE.omegaEV.varphiV}
  \sum_{E \in \mathcal{E}_T} \omega_{TE} \sum_{V \in \mathcal{V}_E} \omega_{EV} \varphi_V
  =\sum_{V \in \mathcal{V}_T} \varphi_V \sum_{E \in \ET,\,V\in\mathcal V_E} \omega_{TE} \omega_{EV}
  = 0,
\end{equation}
where the conclusion is obtained noting that, for each $V\in\VT$, there are two edges $E_1,E_2\in \ET$ such that $V\in\mathcal V_{E_i}$, and that $\omega_{TE_1}\omega_{E_1V}+\omega_{TE_2}\omega_{E_2V}=0$. Identity \eqref{eq:sum.omegaTE.omegaEV.varphiV} translates the well-known fact, for chains (resp.~cochains), that the boundary of a boundary (resp.~coboundary of a coboundary) vanishes: $\partial\circ\partial=0$ (resp.~$\delta\circ\delta=0$).

\subsection{Polynomial spaces}

For any $Y\in\Th\cup\Eh$, we denote by $\Poly{\ell}(Y)$ the space spanned by the restrictions to $Y$ of bivariate polynomials of total degree $\le \ell$, with the convention that $\Poly{\ell}(Y) \coloneq \{0\}$ for $\ell\leq -1$.
We also let $\Poly{0,\ell}(Y)$ denote the subspace of $\Poly{\ell}(Y)$ spanned by polynomials with zero average over $Y$.
For $\ell\in\Zintegers$ and $X\in\Th\cup\Eh$, we denote by $\lproj{\ell}{X}:L^2(X)\to \Poly{\ell}(X)$ the $L^2$-orthogonal projector onto $\Poly{\ell}(X)$. The tensorised version of the above projection, acting on $L^2(X)^2$ or $L^2(X)^{2\times 2}$, is denoted by the bold version $\vlproj{\ell}{X}$ for both spaces. Any potential ambiguity is resolved by the nature of the object being projected.
In what follows, we will also need spaces of piecewise polynomial functions continuous over the boundary of an element or over the mesh skeleton.
Specifically,
for $\bullet\in\{T,h\}$, we denote by $\Poly{\ell}_c(\mathcal{E}_{\bullet})$ the space of continuous functions on $\bigcup_{E \in \mathcal{E}_\bullet} \overline{E}$ whose restriction to any $E\in\mathcal{E}_{\bullet}$ is in $\Poly{\ell}(E)$.
The space of broken polynomials of total degree $\le \ell$ on $\Th$ is denoted by $\Poly{\ell}(\Th)$.

For a smooth enough scalar-valued function $q$, let $\CURL q \coloneq (\partial_2 q,-\partial_1q)^\top$
and notice that,
for all $T\in\Th$, $E\in\ET$ and $r\in \Cspace{1}(\overline{T})$, it holds
\begin{equation}\label{eq:curl.normal}
  (\CURL r)_{|E}\cdot \normal_E=-\partial_{\tangent_E}r_{|E}
\end{equation}
with, as expected, $\partial_{\tangent_E}r_{|E}$ the derivative of $r_{|E}$ in the direction $\tangent_E$.
For every $T\in\Th$, we fix a point $\bvec{x}_T\in T$ such that $T$ contains a ball centered at $\bvec{x}_T$ of diameter uniformly comparable to $h_T$ and, for
any integer $\ell\ge 0$, we define the following subspaces of $\Poly{\ell}(T)^2$:
\[
\Roly{\ell}(T) \coloneq \CURL\Poly{\ell+1}(T), \qquad
\cRoly{\ell}(T) \coloneq (\bvec{x}-\bvec{x}_T)\Poly{\ell-1}(T).
\]
We have the following direct (non-orthogonal) decomposition \cite{Arnold:18}:
\begin{equation}\label{eq:decomp.Pell}
\Poly{\ell}(T)^2 = \Roly{\ell}(T) \oplus \cRoly{\ell}(T).
\end{equation}
The $L^2$-orthogonal projector onto $\Roly{\ell}(T)$ is denoted by $\Rproj{\ell}{T}$.


\section{Discrete BGG diagram}\label{sec:discrete.bgg}

In this section we describe the following discrete counterpart of \eqref{eq:double.complex.cont}, in which $k\ge 0$ is a measure
of the polynomial consistency of the operators:
\begin{equation}\label{eq:double.complex}
  \begin{tikzcd}[column sep=2.5em]
    \text{DS($k$):}
    & 0 \arrow{r}{}
    & \XSgrad{\Th}\arrow{r}{\SGRAD{h}}
    & \XSrot{\Th}\arrow{r}{\SROT{h}}
    & \Poly{k}(\Th)\arrow{r}{}
    & 0 \\
    \text{DDR$(k+1)^2$:}
    & 0 \arrow{r}{}
    & \tXdRgrad{\Th}\arrow{r}{\tGRAD{h}}\arrow[leftrightarrow]{ur}{\Id}
    & \tXdRrot{\Th}\arrow{r}{\tROT{h}}\arrow{ur}{\sskwh}
    & \Poly{k+1}(\Th)^2\arrow{r}{}
    & 0.
  \end{tikzcd}
\end{equation}
The notation for the spaces and operators is inspired by the continuous one \eqref{eq:double.complex.cont}. Specifically, the subscripts ``2'' and ``1'' are used to differentiate $\XSgrad{\Th}$ (the discrete counterpart of $H_2(\Omega)$) from $\XdRgrad{\Th}$ (the discrete counterpart of $H_1(\Omega)$). The same principle is employed to distinguish the gradient acting on $\XSgrad{\Th}$ (denoted by $\SGRAD{h}$) from that acting on $\tXdRgrad{\Th}$ (denoted by $\tGRAD{h}$).
Similarly, the subscripts ``1'' and ``$\VROT$'' identify the discrete curl operators respectively acting on $\tXdRgrad{\Th}$ (denoted by $\SROT{h}$) and $\tXdRrot{\Th}$ (denoted by $\tROT{h}$).
\begin{remark}[Convention for the polynomial degrees]\label{rem:polynomial.degrees}
    The integers $k$ and $k+1$ in the discrete complexes DS($k$) and DDR$(k+1)^2$, and their spaces, correspond to the degree of polynomial consistency of their discrete differential operators. This is expressed, for the discrete rotor in DS($k$), by \eqref{eq:polynomial.consistency.RT} below, and, for DDR$(k+1)^2$, by \cite[Eqs.~(3.13) and (3.20)]{Di-Pietro.Droniou:23}.
    Notice that, according to the BGG construction, the space $\tXdRgrad{\Th}$ acts both as the domain of the gradient operator in the DDR$(k+1)^2$ complex, and as the domain of the scalar rotor operator in the DS($k$) complex.
    In reference to the above convention on the polynomial degrees, this space is of degree $k+1$ in the former interpretation, but only of degree $k$ in the latter.
    See Remark \ref{rem:higher.rot} for a justification of the degree of $\SROT{h}$.
\end{remark}

  \begin{remark}[Three-dimensional construction]
    To develop a three-dimensional counterpart of \eqref{eq:double.complex}, corresponding to a discrete version of the first two rows of the BGG diagram \cite[Eq.~(32)]{Arnold.Hu:21}, one could, for instance, use the three-dimensional DDR complex of \cite{Di-Pietro.Droniou:23} as the bottom sequence and the three-dimensional Stokes complex of \cite{Hanot:23}. ``Vector skew'' and ``trace'' maps will then have to be designed and analysed for this discrete BGG diagram. As done here in the two-dimensional case, the use of some level of serendipity might simplify the construction.
\end{remark}

\subsection{Tensorised discrete de Rham complex}
\label{sec:tensorized.ddr}

The bottom row of \eqref{eq:double.complex} corresponds to the tensorisation of a version of the two-dimensional serendipity discrete de Rham (DDR) complex of \cite{Di-Pietro.Droniou:23*1} of degree $k + 1$.
We briefly recall it hereafter.

\subsubsection{Spaces and interpolators}

We define the following spaces:
\begin{align}\label{eq:def.Xgrad}
  \tXdRgrad{\Th} &\coloneq
  \begin{aligned}[t]
    \Big\{
    &\underline{\bvec{v}}_h=((\bvec{v}_T)_{T\in\Th},(\bvec{v}_E)_{E\in\Eh},(\bvec{v}_V)_{V\in\Vh})\,:
    \\
    &\bvec{v}_T\in\Poly{k-1}(T)^2\quad\forall T\in\Th\,,\quad
    \bvec{v}_E\in\Poly{k}(E)^2\quad\forall E\in\Eh\,,
    \\
    &\bvec{v}_V\in\Real^2\quad\forall V\in\Vh\Big\},
  \end{aligned}
  \\ \label{eq:def.Xrot}
  \tXdRrot{\Th} &\coloneq
  \begin{aligned}[t]
    \Big\{
    &\underline{\bvec{\tau}}_h=((\bvec{\tau}_T)_{T\in\Th},(\bvec{\tau}_E)_{E\in\Eh})\,:
    \\
    &\bvec{\tau}_T\in\Poly{k}(T)^{2\times 2}\quad\forall T\in\Th\,,\quad
    \bvec{\tau}_E\in\Poly{k+1}(E)^2\quad\forall E\in\Eh\Big\}.
  \end{aligned}
\end{align}
In what follows, the restriction of a discrete space to a mesh element or edge $Y\in\Th\cup\Eh$, obtained by collecting the degrees of freedom attached to $Y$ and its boundary, is denoted by replacing $\Th$ with $Y$, e.g.,
\[
\tXdRrot{T} \coloneq
   \Big\{
    \underline{\bvec{\tau}}_h=(\bvec{\tau}_T,(\bvec{\tau}_E)_{E\in\ET})\,:
    \bvec{\tau}_T\in\Poly{k}(T)^{2\times 2}\quad\forall T\in\Th\,,\quad
    \bvec{\tau}_E\in\Poly{k+1}(E)^2\quad\forall E\in\ET\Big\}.
\]
On the other hand, if $\underline{v}_h$ is a vector of degrees of freedom, its restriction to $Y$ is denoted by $\underline{v}_Y$.

\begin{remark}[Use of serendipity]\label{rem:serendipity}
  The DDR complex of degree $k+1$  would involve cell components in $\Poly{k}(T)$ in $\XdRgrad{\Th}$, and in a trimmed space between $\Poly{k}(T)^{2}$ and $\Poly{k+1}(T)^{2}$ in $\XdRrot{\Th}$; see \cite{Di-Pietro.Droniou:23}.
  We instead consider here a serendipity version obtained taking $\eta_T=3$ in \cite[Assumption 12]{Di-Pietro.Droniou:23*1}, which is always possible.
  This choice simplifies the description of the spaces and
  reduces their dimensions while preserving approximation properties.
  See Remark~\ref{rem:comparison} for a discussion of its impact on the Stokes complex of Section \ref{sec:stokes}.
\end{remark}

The interpolators on discrete spaces provide the interpretation of the vector of polynomials representing a smooth enough function. For the spaces defined by \eqref{eq:def.Xgrad} and \eqref{eq:def.Xrot}, the interpolators are $\IdRgrad:\Cspace{0}(\overline{\Omega})^2\to \tXdRgrad{\Th}$ and $\IdRrot:\Cspace{0}(\overline{\Omega})^{2\times 2}\to \tXdRrot{\Th}$  such that
\begin{alignat}{2}
  \label{eq:def.IdRgrad}
  \IdRgrad \bvec{v}\coloneq {}&((\vlproj{k-1}{T}\bvec{v})_{T\in\Th},(\vlproj{k}{E}\bvec{v})_{E\in\Eh},(\bvec{v}(\bvec{x}_V))_{V\in\Vh})&&\qquad\forall \bvec{v}\in \Cspace{0}(\overline{\Omega})^2,\\ \label{eq:def.IdRrot}
  \IdRrot \bvec{\tau}\coloneq {}&((\vlproj{k}{T}\bvec{\tau})_{T\in\Th},(\vlproj{k+1}{E}(\bvec{\tau}\tangent_E))_{E\in\Eh})&&\qquad\forall \bvec{\tau}\in \Cspace{0}(\overline{\Omega})^{2\times 2}.
\end{alignat}
The convention adopted above for restrictions of vectors to mesh elements or edges will also be used for interpolator; for example, $\IdRgradT \bvec{v}=(\vlproj{k-1}{T}\bvec{v},(\vlproj{k}{E}\bvec{v})_{E\in\ET},(\bvec{v}(\bvec{x}_V))_{V\in\VT})$.

\subsubsection{Discrete differential operators}

The description of the tensorised Discrete de Rham complex is completed by the definitions of $\tGRAD{h}$ and $\tROT{h}$, obtained applying suitable reduction and extension maps to the corresponding operators in \cite{Di-Pietro.Droniou:23*1} and further accounting for the projection properties \cite[Eqs.~(6.5) and (6.7)]{Di-Pietro.Droniou:23*1}.
Specifically, for all $\dofh{\bvec v}\in\tXdRgrad{\Th}$, all $E\in\Eh$, and all $T\in\Th$, the discrete edge gradient $\tGE\dof{\bvec v}{E}\in\Poly{k+1}(E)^2$ and discrete element gradient $\tGT\dofT{\bvec v}\in\Poly{k}(T)^{2\times 2}$ are respectively such that
\begin{subequations}\label{eq:def.tGrad}
  \begin{gather}
    \int_E \tGE\underline{\bvec{v}}_E\cdot\bvec{w}
    =-\int_E \bvec{v}_E\cdot \partial_{\tangent_E}\bvec{w}
    + \sum_{V\in\VE}\omega_{EV}\bvec{v}_V\cdot\bvec{w}(\bvec{x}_V)
    \qquad \forall \bvec{w}\in \Poly{k+1}(E)^2, \label{eq:def.tGrad.E}
    \\
    \int_T \tGT\underline{\bvec{v}}_T:\bvec{\zeta}
    = -\int_T \bvec{v}_T\cdot\DIV\bvec{\zeta}
    + \sum_{E\in\ET}\omega_{TE}\int_{E} \bvec{v}_E \cdot (\bvec{\zeta} \normal_E)
    \qquad \forall\bvec{\zeta}\in\Poly{k}(T)^{2\times 2}.\label{eq:def.tGrad.T}
  \end{gather}
  The discrete gradient of $\dofh{\bvec v}$ is then given by
  \begin{equation}\label{eq:def.uGh1}
    \tGRAD{h}\dofh{\bvec v}\coloneq ((\tGT\dof{\boldsymbol{v}}{T})_{T\in\Th},(\tGE\dof{\boldsymbol{v}}{E})_{E\in\Eh})
    \in\tXdRrot{\Th}.
  \end{equation}
\end{subequations}

For all $T\in\Th$, and all $\dofh{\bvec \tau}\in\tXdRrot{\Th}$, the discrete local counterpart of $\bvec{\ROT}$ is $\tROT{T}\dofT{\bvec \tau}\in\Poly{k+1}(T)^2$ such that
\[
\int_T \tROT{T}\underline{\bvec{\tau}}_T\cdot \bvec w
= \int_T \bvec{\tau}_T : \CURL \bvec w
- \sum_{E\in\ET}\omega_{TE}\int_{E} \bvec{\tau}_E \cdot \bvec w
\qquad\forall \bvec w\in\Poly{k+1}(T)^2.
\]
For all $\dofh{\bvec \tau}\in\tXdRrot{\Th}$, we then let $\tROT{h} \dofh{\bvec \tau} \in\Poly{k+1}(\Th)^2 $ be such that
\begin{equation}\label{eq:tROT.h}
  (\tROT{h}\dofh{\bvec \tau})_{|T} \coloneq \tROT{T}\dofT{\bvec \tau}
  \qquad \forall T\in\Th.
\end{equation}

\subsubsection{Skew operator $\sskwh$}

The discrete skew operator $\sskwh:\tXdRrot{\Th}\to\Poly{k}(\Th)$ is obtained applying the skew operator to the element components: For all $\uvec{\tau}_h\in\tXdRrot{\Th}$,
\begin{equation}\label{eq:def.sskwh}
  (\sskwh \underline{\bvec{\tau}}_h)_{|T}
  \coloneq \sskw \bvec{\tau}_T
  \qquad \forall T \in \Th.
\end{equation}


\subsection{Stokes complex}\label{sec:stokes}

\subsubsection{Space and interpolator}

To define the DS($k$) complex in \eqref{eq:double.complex}, we reverse-engineer the space $\XSgrad{\Th}$ to ensure that it contains sufficient information to reconstruct a gradient in $\tXdRgrad{\Th}$.
As explained in, e.g., \cite{Di-Pietro.Droniou.ea:20,Di-Pietro.Droniou:23}, the vertex and edge components of $\tXdRgrad{\Th}$ correspond to the vertex values and $L^2$-orthogonal projections of degree $k$ on edges of functions in $\Poly{k+2}_c(\Eh)^2$. As a consequence, $\XSgrad{\Th}$ should embed enough information to reconstruct, in this space, a tangential and normal gradient on the boundary.
We will see that this is the case for the following choice:
\begin{equation}\label{eq:XSgrad}
  \begin{aligned}
    \XSgrad{\Th}\coloneq \Big\{{}&\underline{q}_h=((q_T)_{T\in\Th},(q_E)_{E\in\Eh},(\Gqen)_{E\in\Eh},(q_V)_{V\in\Vh},(\Gqv)_{V\in\Vh})\,:\\
           {}&q_T\in\Poly{k-2}(T)\quad\forall T\in\Th\,,\quad
           q_E\in\Poly{k-1}(E)\text{ and }\Gqen\in\Poly{k}(E)\quad\forall E\in\Eh\,,\\
           {}&q_V\in\Real\text{ and }\Gqv\in\Real^2\quad\forall V\in\Vh\Big\}.
  \end{aligned}
\end{equation}
We identify the meaning of the polynomial components through the interpolator $\ISgrad:\Cspace{1}(\overline{\Omega})\to \XSgrad{\Th}$ defined by
\begin{equation}\label{eq:def.ISgrad}
  \begin{aligned}
    \ISgrad q\coloneq \Big({}&(\lproj{k-2}{T}q)_{T\in\Th},(\lproj{k-1}{E}q)_{E\in\Eh},(\lproj{k}{E}(\GRAD q\cdot\normal_E))_{E\in\Eh},
    \\
      {}&(q(\bvec{x}_V))_{V\in\Vh},(\GRAD q(\bvec{x}_V))_{V\in\Vh}
      \Big)\qquad\forall q\in \Cspace{1}(\overline{\Omega}).
  \end{aligned}
\end{equation}

  \begin{remark}[Polynomial degrees of the DOFs in $\XSgrad{\Th}$ and $\XSrot{\Th}$]
  The polynomial degrees of the spaces are not related to the $H_1$- or $H_2$-continuity properties mimicked by the spaces. As in finite elements, these continuity properties are mirrored by the order of the derivatives in the \emph{boundary} degrees of freedom.
  In this context, the lower degree of internal DOFs in $\XSgrad{\Th}$ with respect to $\XSrot{\Th}$ is compensated by additional information on the boundary.
  \end{remark}

\begin{remark}[Comparison with the Stokes complex of \cite{Hanot:23}]\label{rem:comparison}
  For the same degree of consistency, the DS complex in \eqref{eq:double.complex} embeds polynomial spaces of one degree lower on each element than the discrete Stokes complex of \cite{Hanot:23}, in the discrete counterparts of $H_2(\Omega)$ and $H_1(\Omega)^2$. This can be explained by the use of serendipity in the present work.
\end{remark}

\subsubsection{Discrete differential operators}

The components of the discrete gradient and rotor for the DS complex are obtained mimicking appropriate integration by parts formulas as described below.

For all $E \in\Eh$ and all $\dof{q}{E}\in\XSgrad{E}$, the discrete tangential gradient $\Gqet \dof{q}{E} \in \Poly{k}(E)$ is such that  \begin{equation}
  \label{eq:def.Gqet}
  \int_E \Gqet\underline{q}_E\,r=-\int_E q_E \partial_{\tangent_E} r + \sum_{V\in\VE}\omega_{EV}\,q_V\,r(\bvec{x}_V)\qquad\forall r\in\Poly{k}(E).
\end{equation}
For all $T\in\Th$ and all $\dof{q}{T}\in\XSgrad{T}$, the discrete element gradient $\nablaT\underline{q}_T\in\Poly{k-1}(T)^2$ satisfies
\begin{equation}
  \label{def:nablaT}
  \int_T\nablaT\underline{q}_T\cdot \bvec w =
  -\int_T q_T\DIV \bvec w
  + \sum_{E\in\ET}\omega_{TE}\int_E q_E\, (\bvec w\cdot\normal_{E})\qquad
  \forall \bvec w\in\Poly{k-1}(T)^2.
\end{equation}
The discrete gradient of $\underline{q}_h\in\XSgrad{\Th}$ is then given by
\begin{equation}\label{eq:def.nablah}
  \SGRAD{h}\dofh{q}=((\nablaT\underline{q}_T)_{T\in\Th},(\Gqet\dof{q}{E}\tangent_E + \Gqen \normal_E))_{E\in\Eh},(\Gqv)_{V\in\Vh})
  \in \XSrot{\Th}.
\end{equation}

\begin{remark}[Notation]
  The notations for discrete gradients associated with $\dof{q}{h}\in\XSgrad{\Th}$ are chosen to distinguish those that are \emph{components} of $\dof{q}{h}$, for which a simple subscript $q$ is used (as in $\Gqen$, $\Gqv$), and those that are \emph{reconstructed} from $\dof{q}{h}$, for which an operator-like notation is used (as in $\nablaT\underline{q}_T$, $\Gqet\dof{q}{E}$).
\end{remark}

\subsubsection{Scalar rot $\SROT{h}$}

For all $\underline{\bvec v}_h\in\tXdRgrad{\Th}$, we define $\SROT{h}\dof{\bvec{v}}{h}$ such that
\begin{equation}\label{eq:SROT.h}
  (\SROT{h}\underline{\bvec v}_h)_{|T}
  = \SROT{T}\dof{\bvec v}{T}
  \qquad \forall T \in \Th,
\end{equation}
where $\SROT{T}\dof{\bvec v}{T} \in\Poly{k}(T)$ satisfies
\begin{equation}\label{eq:def.SROT}
  \int_T\SROT{T}\dof{v}{T}\,r=\int_T \bvec{v}_T\cdot\CURL r - \sum_{E\in\ET}\omega_{TE}\int_{E}(\bvec{v}_E\cdot \tangent_E)\,r\qquad
  \forall r\in\Poly{k}(T).
\end{equation}

  \begin{remark}[Higher-order discrete rotor]\label{rem:higher.rot}
    While the space $\tXdRgrad{\Th}$ contains, in principle, sufficient information to reconstruct a polynomially consistent rotor of degree $k+1$, this reconstruction would fail to ensure the exactness of the tail of the DS($k$) complex.
    In particular, the argument of Lemma \ref{lem:global.exactness} would break down in this case.
  \end{remark}


\section{Homological properties of the Stokes complex}\label{sec:homological_properties}

This section contains lemmas expressing key homological properties of the discrete Stokes complex.

\begin{theorem}[Cohomology of the discrete Stokes complex]\label{th:coho.stokes}
  For all $k \ge 0$, the DS($k$) sequence in \eqref{eq:double.complex} is a complex and its cohomology is isomorphic to the continuous de Rham cohomology.
\end{theorem}

\begin{proof}
  See Section \ref{sec:proof.th.coho}.
\end{proof}

\begin{lemma}[Local commutation properties]\label{lem:loc_com_prop}
  It holds, for all $T\in\Th$,
  \begin{alignat}{2}
    \label{eq:commutation.grad}
    \SGRAD{T} \ISgradT q {}&=\IdRgradT (\GRAD q)  &&\qquad\forall q \in \Cspace{1}(\overline{T}),\\
    \label{eq:commutation.rot}
    \SROT{T} \IdRgradT \bvec v {}&= \lproj{k}{T} (\ROT \bvec v)  &&\qquad\forall \bvec v \in H_2(T)^2.
  \end{alignat}
\end{lemma}

The following polynomial consistency property of $\SROT{T}$ follows immediately from Lemma \ref{lem:loc_com_prop}:
\begin{equation}\label{eq:polynomial.consistency.RT}
    \SROT{T} \IdRgradT \bvec v = \ROT \bvec v  \qquad\forall \bvec v \in \Poly{k+1}(T)^2.
\end{equation}

\begin{remark}[Cochain map property for the interpolators]
  A consequence of Lemma~\ref{lem:loc_com_prop} is the cochain map property for the interpolators, expressed by the commutativity of the following diagram:
  \[
  \begin{tikzcd}[row sep=large]
    H_3(\Omega)\arrow{r}{\GRAD}\arrow{d}{\ISgrad}
    & H_2(\Omega)^2\arrow{r}{\ROT}\arrow{d}{\IdRgrad}
    & H_1(\Omega)\arrow{d}{\lproj{k}{h}}
    \\
    \XSgrad{\Th}\arrow{r}{\SGRAD{h}}
    & \XSrot{\Th}\arrow{r}{\SROT{h}}
    & \Poly{k}(\Th).
  \end{tikzcd}
  \]
  Notice that the continuous Stokes complex in the top row has increased regularity with respect to the one in~\eqref{eq:double.complex.cont} to accommodate the continuity requirements of the interpolators.
\end{remark}

\begin{proof}[Proof of Lemma~\ref{lem:loc_com_prop}]
  i) \emph{Proof of \eqref{eq:commutation.grad}}. Let $q \in \Cspace{1}(\overline{T})$. The equality of the vertex components in \eqref{eq:commutation.grad} is an immediate consequence of the definitions of the interpolators \eqref{eq:def.IdRgrad}, \eqref{eq:def.ISgrad}, and of the discrete gradient \eqref{eq:def.nablah}.
  For all $E\in\ET$ and $r\in\Poly{k}(E)$, we have
  \begin{align*}
    \int_E(\SGRAD{T}\ISgradT q)_E\cdot \tangent_E\,r
    \overset{\eqref{eq:def.nablah}}&=\int_E\Gqet\ISgradT q \, r\\
    \overset{\eqref{eq:def.Gqet},\eqref{eq:def.ISgrad}}&=-\int_E \cancel{\lproj{k-1}{E} q} \, \partial_{\tangent_E}r + \sum_{V\in\VE}\omega_{EV}\,q(\bvec{x}_V)\,r(\bvec{x}_V)
    \overset{\text{IBP}}=\int_E (\GRAD q\cdot\tangent_E) \, r,
  \end{align*}
  which, together with $(\SGRAD{T}\ISgradT q)_E\cdot \normal_E
  \overset{\eqref{eq:def.nablah},\eqref{eq:def.ISgrad}}=\lproj{k}{E} (\GRAD q \cdot\normal_E)$, yields
  \[
  (\SGRAD{T}\ISgradT q)_E=\vlproj{k}{E}(\GRAD q)\overset{\eqref{eq:def.IdRgrad}}=(\IdRgradT \GRAD q)_E,
  \]
  expressing the equality of the edge components in \eqref{eq:commutation.grad}.
  Moving to the element component, for all $T\in\Th$, we have $(\SGRAD{T}\ISgradT q)_T=\nablaT(\ISgradT q)$ and $(\IdRgradT \GRAD q)_T=\vlproj{k-1}{T}(\GRAD q)$.
  Moreover, for all $\bvec w\in\Poly{k-1}(T)^2$,
  \[
  \int_T \nablaT (\ISgradT q) \cdot \bvec w
  \overset{\eqref{def:nablaT},\eqref{eq:def.ISgrad}}=
  -\int_T \cancel{\lproj{k-2}{T}} q \, \DIV \bvec w
  + \sum_{E\in\ET}\omega_{TE}\int_E \cancel{\lproj{k-1}{E}} q\, (\bvec w\cdot\normal_{E})\overset{\text{IBP}}= \int_T \GRAD q \cdot \bvec w,
  \]
  where the removal of the projectors is justified by the respective definitions after observing that $\DIV\bvec{w}\in\Poly{k-2}(T)$ and $\bvec{w}\cdot\normal_{E}\in\Poly{k-1}(E)$ for all $E\in\ET$.
  This relation gives the equality of the element components in \eqref{eq:commutation.grad}, which concludes the proof of this relation.
  \medskip
  \\
  ii) \emph{Proof of \eqref{eq:commutation.rot}}. If $\bvec v\in H_2(T)^2$, we simply observe that, for all $r\in\Poly{k}(T)$,
  \[
  \int_T \SROT{T}(\IdRgradT \bvec v)\, r
  \overset{\eqref{eq:def.SROT},\eqref{eq:def.IdRgrad}}= \int_T \cancel{\vlproj{k-1}{T}} \bvec v \, \CURL r - \sum_{E\in\ET}\omega_{TE}\int_{E}\cancel{\vlproj{k}{E}}\bvec v \cdot \tangent_E\,r
  \overset{\text{IBP}}= \int_T \ROT \bvec v \, r. \qedhere
  \]
\end{proof}

\subsection{Preliminary results for the cohomology analysis}

We state and prove in this section some results
that are required in the proof of Theorem \ref{th:coho.stokes}.

\begin{lemma}[Complex property]
  For all $k \ge 0$, the DS($k$) sequence in \eqref{eq:double.complex} is a complex, i.e.,
  \label{lem:complex.property}
  \begin{align}
    \label{eq:complx.prty.I.grad}
    \ISgrad(\Real){}&\subset\Ker(\SGRAD{h}),\\
    \label{eq:complx.prty.grad.rot}
    \Image(\SGRAD{h}){}&\subset\Ker(\SROT{h}).
  \end{align}
\end{lemma}

\begin{proof}
  i) \emph{Equation \eqref{eq:complx.prty.I.grad}}.
  Straightforward consequence of \eqref{eq:commutation.grad}.
  \medskip
  \\
  ii) \emph{Equation \eqref{eq:complx.prty.grad.rot}}. Let $\dofh{q}\in\XSgrad{\Th}$. For all $r\in\Poly{k}(T)$,
  \begin{align}
    \int_T \nablaT \dof{q}{T} \cdot \CURL r \overset{\eqref{def:nablaT}}&= - \int_T q_T \, \cancel{\DIV\CURL} r + \sum_{E\in\ET}\omega_{TE}\int_E q_E\,(\CURL r \cdot\normal_{E})\nonumber\\
    \overset{\eqref{eq:curl.normal}}& =-\sum_{E\in\ET}\omega_{TE}\int_E q_E\, \partial_{\tangent_E} r_{|E}\nonumber\\
    \overset{\eqref{eq:def.Gqet}}&=\sum_{E\in\ET}\omega_{TE} \int_E \Gqet\underline{q}_E\,r  - \cancel{\sum_{E\in\ET}\omega_{TE}\sum_{V\in\VE}\omega_{EV} q_V \, r(\bvec{x}_V)},
    \label{eq:RTGT=0}
  \end{align}
  where the cancellation in the conclusion comes from \eqref{eq:sum.omegaTE.omegaEV.varphiV} with $(\varphi_V)_{V \in \VE} = (q_V \, r(\bvec{x}_V))_{V \in \VE}$.
  The relation $\SROT{T}\SGRAD{T}\dofT{q}=0$ follows by plugging \eqref{eq:RTGT=0} into the definition \eqref{eq:def.SROT} of $\SROT{T}$ with $\uvec{v}_T= \SGRAD{T}\dofT{q}$.
\end{proof}

\begin{lemma}[Local exactness]
  \label{lem:local.exactness}
  For all $T\in \Th$,
  \begin{align}
    \label{eq:loc.ex.I.grad}
    \ISgradT(\Real){}&=\Ker(\SGRAD{T}),\\
    \label{eq:loc.ex.grad.rot}
    \Image(\SGRAD{T}){}&=\Ker(\SROT{T}).
  \end{align}
\end{lemma}

\begin{proof}
  Let $T\in\Th$. First, notice that Lemma \ref{lem:complex.property} gives the inclusions $\ISgradT(\Real)\subset\Ker(\SGRAD{T})$ and $\Image(\SGRAD{T})\subset\Ker(\SROT{T})$. Hence, only the converse inclusions remain to be proved.\\
  ~\\
  i) \emph{Proof of \eqref{eq:loc.ex.I.grad}}.  Let $\dof{q}{T}\in \XSgrad{T}$ be such that $\SGRAD{T}\dof{q}{T}=\underline{\bvec{0}} \in \XSrot{T}$.
  The definition \eqref{eq:def.nablah} then shows that $\Gqv=\bvec{0}$, $\Gqen=0$, and $\Gqet \dof{q}{E}=0$ for all $V\in\VT$ and $E\in\ET$.
  For $E\in\ET$, take $r=1$ in \eqref{eq:def.Gqet} with $\Gqet \dof{q}{E}=0$ to see that the vertex values of $\dof{q}{T}$ at both edges of $E$ coincide.
  Since $\partial T$ is connected, this gives the equality of all vertex values on $\VT$. Denoting by $C$ their common value, we then go back to \eqref{eq:def.Gqet} with $r$ generic and see that $\int_E q_E \, \partial_{\tangent_E} r=\int_E C \, \partial_{\tangent_E} r$ and, thus, that $q_E=\lproj{k-1}{E}C$ since $\partial_{\tangent_E} r$ spans $\Poly{k-1}(E)$ when $r$ spans $\Poly{k}(E)$.
  Finally, recalling that $\nablaT\dof{q}{T}=\bvec{0}$, it follows from \eqref{def:nablaT} and the fact that $q_E=\lproj{k-1}{E}C$ for all $E\in\ET$, that
  \[
  \int_T q_T\DIV \bvec w
  = \sum_{E\in \ET} \omega_{TE}\int_{E} C \, (\bvec w\cdot\normal_{E})
  \overset{\text{IBP}}=\int_T C\DIV \bvec w\qquad\forall \bvec w\in \Poly{k-1}(T)^2.
  \]
  Since $\DIV : \cRoly{k-1}(T) \rightarrow \Poly{k-2}(T)$ is an isomorphism, this shows that $q_T=\lproj{k-2}{T}C$, hence $\dof{q}{T} = \ISgradT(C)$.%
  \medskip\\
  ii) \emph{Proof of \eqref{eq:loc.ex.grad.rot}}. Let $\dofh{\bvec v}\in\XSrot{T}$ be such that $\SROT{T}\dof{\bvec v}{T}=0$.
  Plugging this condition into \eqref{eq:def.SROT} written for $r=1$, we infer
  \begin{equation*}
    \sum_{E\in\ET}\omega_{TE}\int_{E} \bvec{v}_{E}\cdot\tangent_{E} =0.
  \end{equation*}
  Reasoning as in \cite[Proposition 4.2]{Di-Pietro.Droniou.ea:20} gives $\varphi \in \Poly{k+1}_c(\partial T)$ such that, for all $E\in\ET$, $\partial_{\tangent_E} (\varphi_{|E})=\bvec{v}_E\cdot \tangent_E$.
  We then define $\dof{q}{T}\in\XSgrad{T}$ the following way:
  \begin{subequations}\label{eq:loc.exact.def.qT}
    \begin{alignat}{2}
      \label{eq:def.qT.in.proof.V}
      q_V&=\varphi(\bvec{x}_V) &\qquad& \forall V\in\VT, \\
      \label{eq:def.qT.in.proof.GqV}
      \Gqv&=\bvec{v}_V &\qquad& \forall V\in\VT, \\
      \label{eq:def.qT.in.proof.E}
      q_E&=\lproj{k-1}{E}\varphi &\qquad& \forall E\in\ET, \\
      \label{eq:def.qT.in.proof.Gqen}
      \Gqen&=\bvec{v}_E\cdot\normal_E &\qquad& \forall E\in\ET,
    \end{alignat}
    and $q_T\in\Poly{k-2}(T)$ is such that
    \begin{equation} \label{eq:def.qT.in.proof.T}
      \int_T q_T \DIV \bvec w
      = -\int_T \bvec{v}_T \cdot \bvec w
      + \sum_{E\in\ET}\omega_{TE}\int_E q_E\, (\bvec w\cdot\normal_{E})
      \qquad \forall \bvec w \in\cRoly{k-1}(T).
    \end{equation}
  \end{subequations}
  Let us show that $\SGRAD{T}\dof{q}{T} = \dof{\bvec v}{T}$. The equality at the vertices and of the normal components on the edges is an immediate consequence of \eqref{eq:def.qT.in.proof.GqV} and \eqref{eq:def.qT.in.proof.Gqen} together with the definition \eqref{eq:def.nablah} of $\SGRAD{h}$.
  To prove that, for all $E\in\ET$, $\Gqet \dof{q}{E}=\bvec{v}_E\cdot\tangent_E$ we write, for $r\in\Poly{k}(E)$
  \begin{align*}
    \int_E \Gqet \dof{q}{E}\, r
    \overset{\eqref{eq:def.Gqet}}&=  -\int_E q_E\, \partial_{\tangent_E} r
    + \sum_{V\in\VE}\omega_{EV} \,q_V\, r(\bvec{x}_V)
    \\
    \overset{\eqref{eq:def.qT.in.proof.E},\eqref{eq:def.qT.in.proof.V}}&=
    -\int_E \cancel{\lproj{k-1}{E}}\varphi_{|E}\,\partial_{\tangent_E}  r
    + \sum_{V\in\VE} \omega_{EV} \varphi(\bvec{x}_V) r(\bvec{x}_V)
    \\
    \overset{\text{IBP}}&= \int_E \partial_{\tangent_E} (\varphi_{|E})\, r
    =\int_E (\bvec{v}_E\cdot\tangent_E) \,r,
  \end{align*}
  showing that
  \begin{equation}\label{eq:equality.tang.comp}
    \Gqet \dof{q}{E} = \bvec{v}_E\cdot\tangent_E
  \end{equation}
  and thus that $\bvec v_E = (\SGRAD{T}\dofT{q})_E$.
  Finally, we consider the components on $T$.
  For $\bvec{z} \in \Roly{k-1}(T)$, there exists $r\in\Poly{k}(T)$ such that $\bvec{z}=\CURL r$. Recalling that $\SROT{T}\dof{\bvec v}{T}=0$, we infer that
  \begin{equation}\label{eq:proof.loc.ex.un}
    \begin{aligned}
      \int_T \bvec{v}_T \cdot \bvec{z} \overset{\eqref{eq:def.SROT}}&= \sum_{E\in\ET}\omega_{TE}\int_{E}  (\bvec{v}_{E}\cdot\tangent_{E}) \, r
      \\
      \overset{\eqref{eq:equality.tang.comp}}&=
      \sum_{E\in\ET}\omega_{TE}\int_{E} \Gqet \dof{q}{E} \, r
      \\
      \overset{\eqref{eq:def.Gqet},\, \eqref{eq:sum.omegaTE.omegaEV.varphiV}}&=-\sum_{E\in\ET}\omega_{TE} \int_E q_E \partial_{\tangent_E} r_{|E} + \cancel{\sum_{E\in\ET}\omega_{TE} \sum_{V\in\VE}\omega_{EV}\,q_V\,r(\bvec{x}_V)}
      \\
      \overset{\eqref{eq:curl.normal}}&=\sum_{E\in\ET}\omega_{TE} \int_E q_E \, (\bvec{z}\cdot \normal_{E}) - \int_T q_T \DIV \bvec{z},
    \end{aligned}
  \end{equation}
  where, the last equality, the introduction of $\int_T q_T \DIV \bvec{z}$ is justified by $\DIV\bvec{z} =\DIV\CURL r = 0$.
  Adding together  \eqref{eq:proof.loc.ex.un} and \eqref{eq:def.qT.in.proof.T} and recalling that $\bvec z +\bvec w$ spans $\Poly{k-1}(T)^2$ as $(\bvec w,\bvec z)$ spans $\Roly{k-1}(T)\times\cRoly{k-1}(T)$ (see \eqref{eq:decomp.Pell}), we infer that $\dof{q}{T}$ satisfies \eqref{eq:def.qT.in.proof.T} for all $\bvec w\in\Poly{k-1}(T)^2$.
  The definition \eqref{def:nablaT} of $\nablaT\dof{q}{T}$ then yields $\bvec{v}_T=\nablaT\dof{q}{T}=(\SGRAD{T}\dof{q}{T})_T$, which concludes the proof.
\end{proof}

\begin{lemma}[Global exactness of the tail of DS($k$)]
  \label{lem:global.exactness}
  It holds
  \begin{align}
    \label{eq:loc.ex.rot}
    \Image(\SROT{h}) {}&= \Poly{k}(\Th).
  \end{align}
\end{lemma}

\begin{proof}
  Let $z_h\in\Poly{k}(\Th)$. By the surjectivity of $\ROT : H_1(\Omega)^2 \to L^2(\Omega)$ (which comes from the surjectivity of $\DIV:H_1(\Omega)^2 \to L^2(\Omega)$, see \cite[Lemma 8.3]{Di-Pietro.Droniou:20} and references therein, together with the fact that $\ROT$ is the divergence of the rotated vector field),
  there exists $\bvec w\in H_1(\Omega)^2$ such that $\ROT \bvec w = z_h$.
  We then define $\dofh{\bvec v}\in\XSrot{\Th}$ such that
  \begin{equation}\label{eq:def.vT.global.exactness}
    \bvec{v}_V= \bvec 0 \qquad\forall V\in\Vh,\qquad \bvec{v}_E= \vlproj{k}{E}\bvec w\quad\forall E\in\Eh,\qquad
    \bvec{v}_T=\vlproj{k-1}{T}\bvec w \quad\forall T\in\Th.
  \end{equation}
  Then, for all $T\in\Th$ and all $r\in\Poly{k}(T)$,
  \begin{align*}
    \int_T\SROT{h}\underline{v}_h\,r\overset{\eqref{eq:def.SROT}}
    &=\int_T \bvec{v}_T\cdot\CURL r - \sum_{E\in\ET}\omega_{TE}\int_{E}(\bvec{v}_E\cdot \tangent_E)\,r\\
    \overset{\eqref{eq:def.vT.global.exactness}}&=
    \int_T \cancel{\vlproj{k-1}{T}} \bvec w\cdot\CURL r - \sum_{E\in\ET}\omega_{TE}\int_{E} \cancel{\lproj{k}{E}} (\bvec w\cdot \tangent_E)\,r
    \overset{\text{IBP}} =\int_T \ROT \bvec w \,r = \int_T (z_h)_{|T}\, r.
  \end{align*}
  Thus, $\SROT{T}\dof{\bvec v}{T} = (z_h)_{|T}$ for all $T\in\Th$, so $\SROT{h}\dofh{\bvec v}=z_h$.
\end{proof}

\subsection{Proof of Theorem \ref{th:coho.stokes}}
\label{sec:proof.th.coho}

In this section we prove that the cohomology of the DS($k$) complex is isomorphic to the cohomology of DDR(0), the (scalar) DDR complex of degree $0$ (see \cite[Section 4.1]{Di-Pietro.Droniou.ea:23}), which is in turn isomorphic to that of the continuous de Rham complex \cite[Lemma 4]{Di-Pietro.Droniou.ea:23}.
We use the framework of \cite{Di-Pietro.Droniou:23*1} and create so-called extension and reduction cochain maps satisfying \cite[Assumption 1]{Di-Pietro.Droniou:23*1}:
\begin{equation}\label{eq:doublecomplex_ddr0_stokes}
  \begin{tikzcd}[column sep=3.5em]
    \text{DS($k$):}\hspace*{-3em}
    & 0 \arrow{r}{\ISgrad} & \XSgrad{\Th}\arrow{r}{\SGRAD{h}} \arrow[d, bend left, "\Rgradh"] & \XSrot{\Th}\arrow{r}{\SROT{h}}\arrow[d, bend left, "\Rroth"]&\Poly{k}(\Th)\arrow{r}{} \arrow[d, bend left, "\lproj{0}{h}"] & 0\\
    \text{DDR(0):}\hspace*{-3em}
    & 0 \arrow{r}{\Igrad{0}{h}}  & \Xgrad{0}{h}\arrow{r}{\uGh{0}} \arrow[u, bend left, "\Egradh"] & \Xcurl{0}{h}\arrow[u, bend left, "\Eroth"]\arrow{r}{\uCh{0}}&\Poly{0}(\Th)\arrow{r}{} \arrow[u, bend left, "i"]& 0.
  \end{tikzcd}
\end{equation}
We recall that the discrete $H_1(\Omega)$ and $\boldsymbol{H}_{\CURL}(\Omega)$ spaces in the DDR(0) are respectively given by
\[
\Xgrad{0}{h} \coloneq \left\{(q_V)_{V\in\Vh}\,:\,q_V\in\Real\quad\forall V\in\Vh\right\},\quad
\Xcurl{0}{h} \coloneq \left\{(v_E)_{E\in\Eh}\,:\,v_E\in\Real\quad\forall E\in\Eh\right\},
\]
with discrete gradient and curl operators respectively defined as
\begin{gather}\label{eq:def.uGh0}
  \uGh{0}\dof{q}{h} \coloneq \left(\GE{0}\dof{q}{E}=\frac{1}{h_E}\sum_{V\in\VE}\omega_{EV}q_V\right)_{E\in\Eh}\quad\forall \dof{q}{h}\in\Xgrad{0}{h},
  \\ \label{eq:def.uCh0}
  \uCh{0}\dof{v}{h} \coloneq \left(-\frac{1}{|T|}\sum_{E\in\ET}\omega_{TE}h_E v_E\right)_{T\in\Th}\quad\forall \dof{v}{h}\in\Xcurl{0}{h}
\end{gather}
and interpolator $\Igrad{0}{h} : \Cspace{0}(\overline{\Omega}) \to \Xgrad{0}{h}$ such that
\[
\Igrad{0}{h}q \coloneq (q(\bvec{x}_V))_{V\in\Vh}
\qquad \forall q \in \Cspace{0}(\overline{\Omega}).
\]

The reduction maps in \eqref{eq:doublecomplex_ddr0_stokes} are such that
\begin{alignat}{2}
  \label{eq:def.Rgrad}
  \Rgradh \dofh{q} &= (q_V)_{V \in \Vh} && \qquad \forall \dofh{q}\in\XSgrad{\Th}, \\
  \label{eq:def.Rrot}
  \Rroth  \dofh{\bvec{v}} &= (\lproj{0}{E} (\bvec{v}_E\cdot \tangent_E))_{E \in \Eh} &&\qquad \forall \dofh{\bvec{v}} \in \XSrot{\Th}.
\end{alignat}
Let us describe the extension maps. $\Egradh$ is such that, for all $\dofh{q}\in\Xgrad{0}{h}$,
\begin{subequations}\label{eq:def.Egradh}
  \begin{equation}
    \label{eq:def.extensiongrad}
    \Egradh \dofh{q} \coloneq  ((\Egrad{k-2}{T}\dof{q}{T})_{T\in\Th},(\Egrad{k-1}{E}\dof{q}{E})_{E\in\Eh},(0)_{E\in\Eh},(q_V)_{V\in\Vh},(\bvec 0)_{V \in\Vh}),
  \end{equation}
  where, $\Egrad{k-1}{E}\dof{q}{E}\in \Poly{k-1}(E)$ and $\Egrad{k-2}{T} \dof{q}{T} \in \Poly{k-2}(T)$ are respectively defined by
  \begin{gather}\label{eq:def.EgradE}
    \int_E \Egrad{k-1}{E}\dof{q}{E} \partial_{\tangent_E} r
    = -\int_E \GE{0}\dof{q}{E}r
    + \sum_{V\in\VE} \omega_{EV} \, q_V \, r(\bvec{x}_V)
    \qquad \forall r\in\Poly{k}(E),
    \\ \label{eq:def.EgradT}
    \int_T \Egrad{k-2}{T} \dof{q}{T} \DIV \bvec w
    = -\int_T \boldsymbol{\gamma}_{T}^0\uGT{0}\dof{q}{T} \cdot \bvec w
    + \sum_{E\in\ET} \omega_{TE} \int_E  \Egrad{k-1}{E}\dof{q}{E} (\bvec w \cdot \boldsymbol n_{E})
    \qquad \forall \bvec w\in \boldsymbol{\mathcal{R}}^{c,k-1}(T),
  \end{gather}
\end{subequations}
with $\boldsymbol{\gamma}_{T}^0$ given \cite[Eq.~(11)]{Di-Pietro.Droniou:21} with $k=0$.
We notice that $\Egrad{k-1}{E}\dof{q}{E}$ is well-defined because the right-hand side of \eqref{eq:def.EgradE} vanishes when $\partial_{\tangent_E} r=0$, and that it only depends on $\dof{q}{E}$, and is therefore fully known when used to define $\Egrad{k-2}{T} \dof{q}{T}$.

The extension $\Eroth$ is such that, for all $\dofh{\bvec v}\in\Xcurl{0}{h}$,
\begin{subequations}\label{eq:def.Eroth}
  \begin{equation}
    \label{eq:def.extensionrot}
    \Eroth\dofh{\bvec v} \coloneq  ((\Erot{k-1}{T}\dof{\bvec v}{T})_{T\in\Th},(v_E\tangent_E)_{E\in\Eh},(\bvec{0})_{V\in\Vh}),
  \end{equation}
  where, for all $T\in \Th$, $\Erot{k-1}{T}\dof{\bvec v}{T}\in \Poly{k-1}(T)^2$ is defined by
  \begin{multline}\label{eq:def.ErotT}
    \int_T \Erot{k-1}{T}\dof{\bvec v}{T}\cdot (\CURL r+ \bvec w) = \int_{T} C^0_T \dof{\bvec v}{T}\, r +\sum_{E\in\ET} \omega_{TE} \int_E v_E \, r +  \int_T \boldsymbol{\gamma}_{T}^0\dof{\bvec v}{T}\cdot \bvec w \\
    \forall (r,\bvec w) \in\Poly{k}(T)\times\cRoly{k-1}(T).
  \end{multline}
\end{subequations}
Notice that, for $ \Erot{k-1}{T}\dof{\bvec v}{T}$ to be well defined, we must ensure that the right-hand side of \eqref{eq:def.ErotT} vanishes when applied to $\bvec w=\bvec 0$ and $r$ such that $\CURL r=\bvec 0$. This holds true by \eqref{eq:def.uCh0} since $\Ker \CURL = \Poly{0}(T)$.

\begin{remark}[Design of extension operators]
  The approach to defining extension operators is well known, as similar operators have already been introduced in \cite{Di-Pietro.Droniou.ea:23,Di-Pietro.Droniou:23*1,Bonaldi.Di-Pietro.ea:25}. The key idea behind these definitions is that applying a high-order discrete calculus operator to the extension of a vector of low-order polynomials should yield the lowest-order discrete calculus operator applied to this vector, ensuring the cochain map property. This principle is illustrated in \eqref{eq:cochainext2}, which results from the construction of $\Egrad{k-1}{E}\dof{q}{E}$ (see \eqref{eq:cochain.extension.Egrad}).
\end{remark}

\begin{lemma}[Cochain property]\label{lem:cochainstokes}
  The extensions and reductions are cochain maps, that is:
  \begin{alignat}{2}
    \label{eq:cochainred2}
    \uGh{0}\Rgradh\dofh{q} &= \Rroth\SGRAD{h}\dofh{q} &&\qquad\forall \dofh{q}\in\XSgrad{\Th},\\
    \label{eq:cochainred3}
    C^{0}_h\Rroth\dofh{\bvec v} &= \lproj{0}{h}\SROT{h}\dofh{\bvec v} &&\qquad\forall \dofh{\bvec v}\in \XSrot{\Th},\\
    \label{eq:cochainext2}
    \SGRAD{h}\Egradh\dofh{q} &= \Eroth\uGh{0}\dofh{q} &&\qquad\forall \dofh{q}\in \Xgrad{0}{h},\\
    \label{eq:cochainext3}
    \SROT{h}\Eroth\dofh{\bvec v} &= C^{0}_h\dofh{\bvec v} &&\qquad\forall \dofh{\bvec v}\in\Xcurl{0}{h}.
  \end{alignat}
\end{lemma}

\begin{proof}
  i) \emph{Cochain map property for the reductions.}
  To prove \eqref{eq:cochainred2}, let $\dofh{q}\in\XSgrad{\Th}$. We have
  \[
  \uGh{0}\Rgradh\dofh{q}
  \overset{\eqref{eq:def.Rgrad}}= \uGh{0}((q_V)_{V\in\Vh})
  \overset{\eqref{eq:def.uGh0}} = \left(\frac{1}{h_E}\sum_{V\in\VE}\omega_{EV}\,q_V\right)_{E\in\Eh},
  \]
  and \eqref{eq:cochainred2} therefore follows by noticing that, for all $E\in\Eh$,
  \begin{equation*}
    (\Rroth\SGRAD{h}\dofh{q})_E
    \overset{\eqref{eq:def.Rrot},\,\eqref{eq:def.nablah}}=
    \lproj{0}{E} \Gqet \dof{q}{E}
    = \frac{1}{h_E}\int_E\Gqet \dof{q}{E}
    \overset{\eqref{eq:def.Gqet} \text{ with $r=1$}}=
    \frac{1}{h_E}\sum_{V\in\VE}\omega_{EV}\,q_V.
  \end{equation*}
  We now prove \eqref{eq:cochainred3}. Let $\dofh{\bvec v}\in\XSrot{\Th}$, take $T\in\Th$, and notice first that
  \begin{align*}
    (\uCh{0}\Rroth\dofh{\bvec v})_T
    \overset{\eqref{eq:def.Rrot}}= C_T^0(( \vlproj{0}{E} \bvec{v}_E \cdot\tangent_E)_{E\in\Eh})
    \overset{\eqref{eq:def.uCh0}}={}& -\frac{1}{|T|}\sum_{E \in \ET} \omega_{TE}h_E \lproj{0}{E}(\bvec{v}_E \cdot\tangent_E)\\
    ={}& -\frac{1}{|T|}\sum_{E \in \ET} \omega_{TE}\int_E \bvec{v}_E \cdot\tangent_E.
  \end{align*}
  We then obtain \eqref{eq:cochainred3} by writing
  \begin{equation*}
    \lproj{0}{T}\SROT{T}\dofT{\bvec v}
    = \frac{1}{|T|}\int_T \SROT{T}\dofT{\bvec v}
    \overset{\eqref{eq:def.SROT}\text{ with $r=1$}}= - \frac{1}{|T|}\sum_{E\in\ET}\omega_{TE}\int_{E}\bvec v_E\cdot \tangent_E.
  \end{equation*}
  \\
  ii) \emph{Cochain map property for the extensions}.
  To prove \eqref{eq:cochainext2}, let $\dofh{q}\in\Xgrad{0}{h}$. For all $V\in \Vh$,
  \begin{equation*}
    (\SGRAD{h}\Egradh\dofh{q})_V
    \overset{\eqref{eq:def.nablah},\,\eqref{eq:def.extensiongrad}}=\bvec 0
    \overset{\eqref{eq:def.extensionrot}}= (\Eroth\uGh{0}\dofh{q})_V.
  \end{equation*}
  For all $E\in\Eh$, the normal components of $(\SGRAD{h}\Egradh\dofh{q})_E$ and $(\Eroth\uGh{0}\dofh{q})_E$ both vanish by \eqref{eq:def.extensiongrad}-\eqref{eq:def.nablah} and \eqref{eq:def.extensionrot}, respectively.
  Let us show the equality of respective tangential components. We have $(\Eroth\uGh{0}\dofh{q})_E= \GE{0}\dof{q}{E}\tangent_E$ and, for all $r\in\Poly{k}(E)$,
  \begin{align}
    \int_E (\SGRAD{h}\Egradh\dofh{q})_E\cdot\tangent_E\, r \overset{\eqref{eq:def.nablah}} &=\int_E \Gqet \EgradE \dof{q}{E}\, r \nonumber\\
    \overset{\eqref{eq:def.Gqet}}&= -\int_E \Egrad{k-1}{E}\dof{q}{E} \,\partial_{\tangent_E} r + \sum_{V\in\VE} \omega_{EV} \, q_V \, r(\bvec{x}_V)\nonumber\\
    \overset{\eqref{eq:def.EgradE}}&= \int_E \GE{0}\dof{q}{E} \,r
    =\int_E (\Eroth\uGh{0}\dofh{q})_E\cdot \tangent_E\,r.
  \label{eq:cochain.extension.Egrad}
  \end{align}
  This proves the equality of the edge components in \eqref{eq:cochainext2}.
  Consider now $T\in\Th$ and let us prove that the element components in \eqref{eq:cochainext2} coincide. For all $(\bvec v, \bvec w)\in\boldsymbol{\mathcal{R}}^{k-1}(T)\times \boldsymbol{\mathcal{R}}^{c,k-1}(T),$ letting $r\in\Poly{k}(T)$ be such that $\bvec v=\CURL r$, we have, on one hand,
  \[
  \begin{aligned}
    &\int_T \nablaT\EgradT\underline{q}_T \cdot (\CURL r+\bvec w)
    \\
    &\qquad
    \begin{aligned}[t]
      \overset{\eqref{def:nablaT}}&= -\int_T \Egrad{k-2}{T}\dof{q}{T}\DIV( \cancel{\CURL r} + \bvec w)+ \sum_{E\in\ET} \omega_{TE} \int_E  \Egrad{k-1}{E}\dof{q}{E} (\CURL r +\bvec w)  \cdot \boldsymbol n_{E}
      \\
      \overset{\eqref{eq:def.EgradT},\eqref{eq:curl.normal}}&=
      \int_T \boldsymbol{\gamma}_{T}^0\uGT{0}\dof{q}{T} \cdot \bvec w -\sum_{E\in\ET} \omega_{TE} \int_E \Egrad{k-1}{E}\dof{q}{E}\,\partial_{\tangent_E} r
      \\
      \overset{\eqref{eq:def.EgradE}}&=
      \int_T \boldsymbol{\gamma}_{T}^0\uGT{0}\dof{q}{T} \cdot \bvec w+\sum_{E\in\ET} \omega_{TE} \left(\int_E \GE{0}\dof{q}{E}r - \sum_{V\in\VE} \omega_{EV} \, q_V \, r(\bvec{x}_V)\right)
      \\
      \overset{\eqref{eq:sum.omegaTE.omegaEV.varphiV}}&=
      \int_T \boldsymbol{\gamma}_{T}^0\uGT{0}\dof{q}{T} \cdot \bvec w +\sum_{E\in\ET} \omega_{TE} \int_E \GE{0}\dof{q}{E}r.
    \end{aligned}
  \end{aligned}
  \]
  On the other hand,
  \[
  \int_T \Erot{k-1}{T}\uGT{0}\dofT{q}\cdot (\CURL r +\bvec w) \overset{\eqref{eq:def.ErotT}}
  = \int_E \cancel{C^{0}_T \uGT{0}\dofT{q}} \,r +\sum_{E\in\ET} \omega_{TE} \int_E \GE{0}\dof{q}{E} r + \int_T \boldsymbol{\gamma}_{T}^0\uGT{0}\dofT{q}\cdot \bvec w,
  \]
  where the cancellation in the first line is a consequence of the complex property of the DDR sequence.
  The components on $T$ of both sides of \eqref{eq:cochainext2} coincide, which concludes the proof of this relation.

  Finally, to prove \eqref{eq:cochainext3}, we write, for all $T\in\Th$ and all $r\in \Poly{k}(T)$,
  \[
  \begin{aligned}
    \int_T \SROT{T}\ErotT\dofT{\bvec v} \,r \overset{\eqref{eq:def.SROT},\eqref{eq:def.extensionrot}}&=\int_T \Erot{k-1}{T}\dofT{\bvec v} \cdot\CURL r - \sum_{E\in\ET}\omega_{TE}\int_{E}(v_E\tangent_E)\cdot\tangent_E\,r\\
    &=\int_T \Erot{k-1}{T}\dofT{\bvec v} \cdot\CURL r - \sum_{E\in\ET}\omega_{TE}\int_{E}v_E \,r
    \overset{\eqref{eq:def.ErotT}}=\int_T C^{0}_T\dofT{\bvec v} r.\qedhere
  \end{aligned}
  \]
\end{proof}

\begin{lemma}[Pseudo-inversion properties of extensions and reductions]
  \label{lem:C2}
  It holds
  \begin{align}
    \label{eq:C2:1}
    (\Egradh \Rgradh-\Id_{\XSgrad{\Th}})(\Ker(\SGRAD{h})) &= \{ \underline{0} \}, \\
    \label{eq:C2:2}
    (\Eroth\Rroth - \Id_{\XSrot{\Th}})(\Ker(\SROT{h}))&\subset \Image(\SGRAD{h}),\\
    \label{eq:C2:3}
    (\lproj{0}{h} - \Id_{\Poly{k}(\Th)})(\Poly{k}(\Th))&\subset \Image(\SROT{h}).
  \end{align}
\end{lemma}

\begin{proof}
  i) \emph{Proof of \eqref{eq:C2:1}}. The proof is a straightforward  adaptation of \cite[Lemma 8]{Di-Pietro.Droniou.ea:23} restricted to two dimensions, using the local exactness of the DS complex \eqref{eq:loc.ex.I.grad}.
  \medskip\\
  ii) \emph{Proof of \eqref{eq:C2:2}}. Let $\dofh{\bvec v}\in\XSrot{\Th}$ be such that $\SROT{h}\dofh{\bvec v}=0$. By the cochain map property, $\Eroth\Rroth \dofh{\bvec v}\in\Ker\SROT{h}$.
  Let $T\in\Th$. By the local exactness property \eqref{eq:loc.ex.grad.rot}, there exists $\dof{q}{T}\in\XSgrad{T}$ such that
  \begin{equation}
    \label{eq:proof.C2.def.x.T}
    \SGRAD{T}\dof{q}{T}=\ErotT\RrotT \dof{\bvec v}{T}-\dof{\bvec v}{T}.
  \end{equation}
  Taking an arbitrary $V_0\in\VT$ and making the substitution $\dof{q}{T} \gets \dof{q}{T} -\ISgradT q_{V_0}$, the polynomial consistency property \eqref{eq:commutation.grad} shows that \eqref{eq:proof.C2.def.x.T} is still valid. We can therefore assume
  in the following that one of the vertex values $q_{V_0}$ of $\dof{q}{T}$ vanishes.
  Let us show that, for all $V\in\VT$ and $E\in\ET$, the components $\Gqv$, $\Gqen$, $q_V$ and $q_E$ of $\dof{q}{T}$ do not depend on $T$.
  For $V\in\VT$, by definition \eqref{eq:def.extensionrot} of $\Eroth$,  \eqref{eq:proof.C2.def.x.T} gives $\Gqv=(\SGRAD{T}\dof{q}{T})_V=-\bvec{v}_V$, which does not depend on $T$.
  Let $E\in\ET$. By \eqref{eq:proof.C2.def.x.T} we have
  \begin{equation}\label{eq:loc.ex.Gedge}
    \Gqet\dof{q}{E}\tangent_E +\Gqen\normal_E=\lproj{0}{E}(\bvec{v}_E\cdot\tangent_E)\tangent_E -\bvec{v}_E.
  \end{equation}
  Taking the dot product with $\normal_E$, we infer that $\Gqen=-\bvec{v}_E\cdot\normal_E$ only depends on $E$. Moreover, taking the dot product of \eqref{eq:loc.ex.Gedge} with $\tangent_E$ and applying $\lproj{0}{E}$, we obtain $\lproj{0}{E}\Gqet\dof{q}{T}=0$, from which we deduce
  \[
  0=\lproj{0}{E}\Gqet\dof{q}{T}\overset{\eqref{eq:def.Gqet} \text{ with $r=1$}}= \frac{1}{h_E}\sum_{V\in\VE} \omega_{EV}\, q_V.
  \]
  Since the relative orientations of the two vertices of $E$ are opposite, this shows that the two vertex values $(q_V)_{V\in\VE}$ are equal. We can make the same observation on each $E\in\ET$ and, since $\partial T$ is connected, we infer that all $(q_V)_{V\in\VT}$ are equal. As we have assumed that at least one of the vertex values of $\dof{q}{T}$ vanishes, this means that all vertex values of this vector vanish, and are therefore independent of $T$.
  Taking the dot product of \eqref{eq:loc.ex.Gedge} with $\tangent_E$ yields $\Gqet\dof{q}{E}=\lproj{0}{E}(\bvec{v}_E\cdot\tangent_E) -\bvec{v}_E\cdot\tangent_E$ and thus, for all $r\in\Poly{k}(E)$, recalling that all vertex values $(q_V)_{V\in\VT}$ vanish,
  \[
  \int_E (\lproj{0}{E}(\bvec{v}_E\cdot\tangent_E) -\bvec{v}_E\cdot\tangent_E)\,r=\int_E \Gqet\dof{q}{E}\,r \overset{\eqref{eq:def.Gqet}}=-\int_E q_E \, \partial_{\tangent_E} r,
  \]
  showing that $q_E$ depends only on $E$. Hence, for all $E\in\ET$, $\dof{q}{E}=(q_E,-\bvec{v}_E\cdot\normal_E,(0)_{V\in\VE},(-\bvec{v}_V)_{V\in\VE})$ does not depend on $T$.
  The fact that the vertex and edge values of $\dof{q}{T}$ do not depend on $T$ allows us to glue all these local vectors into a global one $\dof{q}{h}\in \XSgrad{\Th}$ (avoiding any risk of multiple definitions of vertex/edge values coming from different elements), such that $\SGRAD{h}\dof{q}{h}=\Eroth\Rroth \dof{\bvec v}{h}-\dof{\bvec v}{h}$.
  \medskip\\
  iii) \emph{Proof of \eqref{eq:C2:3}}. This property is an immediate consequence of the global exactness of $\SROT{h}$, see \eqref{eq:loc.ex.rot}.
\end{proof}

\begin{proof}[Proof of Theorem \ref{th:coho.stokes}]
  The conditions (C1), (C2) and (C3) of \cite[Assumption 1]{Di-Pietro.Droniou:23*1} are satisfied. Indeed, (C1) is a straightforward consequence of the definitions of the operators:
  \[
  \begin{alignedat}{2}
    \Rgradh \Egradh \dofh{q}
    \overset{\eqref{eq:def.extensiongrad},\,\eqref{eq:def.Rgrad}}&=
    \dofh{q} &\qquad& \forall \dofh{q}\in\Xgrad{0}{h},
    \\
    \Rroth \Eroth \dofh{\bvec v}
    \overset{\eqref{eq:def.extensionrot},\,\eqref{eq:def.Rrot}}&=
    \dofh{\bvec v} &\qquad& \forall \dofh{\bvec v}\in \Xcurl{0}{h},
  \end{alignedat}
  \]
  (C2) is established in Lemma~\ref{lem:C2}, and the cochain maps property (C3) is proved in Lemma~\ref{lem:cochainstokes}.
  Thus, the isomorphism property between the DS($k$) complex and the DDR(0) complex follows from \cite[Proposition 2]{Di-Pietro.Droniou:23*1}, and the theorem follows from \cite[Lemma 4]{Di-Pietro.Droniou.ea:23}.
\end{proof}


\section{Discrete twisted and BGG complexes}\label{sec:twisted_and_bgg_complexes}

\subsection{Anti-commutation property of $\sskwh$}

\begin{lemma}[Anti-commutativity]\label{lem:anticomm}
  The diagram \eqref{eq:double.complex} is anti-commutative, that is:
  \[
  \sskwh\circ \,\tGRAD{h}=-\SROT{h}.
  \]
\end{lemma}

\begin{proof}
  Let $\underline{\bvec v}_h\in\tXdRgrad{\Th}$.
  Recalling the definitions \eqref{eq:def.sskwh} of $\sskwh$,  \eqref{eq:def.tGrad} of $\tGRAD{h}\underline{\bvec{v}}_h$ and \eqref{eq:def.SROT} of $\SROT{T}\dof{v}{T}$, we have to show that, for any $T\in\Th$,
  \[
  \sskw(\tGT\underline{\bvec{v}}_T)=-\SROT{T}\dof{v}{T}.
  \]
  Take $r\in\Poly{k}(T)$ and set
  \[
  \bvec{\zeta} \coloneq r \begin{pmatrix} 0&1\\-1&0 \end{pmatrix}
  \in\Poly{k}(T)^{2\times 2}.
  \]
  We have $\tGT\dof{v}{T}:\bvec{\zeta} = \sskw(\tGT\dof{v}{T})r$,
  $\bvec \DIV \bvec{\zeta}=(\partial_2 r,-\partial_1 r)^\top=\CURL r$,
  and $\bvec{\zeta}\normal_E = r\tangent_E$.
  Expressing \eqref{eq:def.tGrad.T} with this choice of $\bvec{\zeta}$ therefore yields
  \[
  \int_T \sskw(\tGT\underline{\bvec{v}}_T) \, r
  = -\int_T \bvec{v}_T\cdot \CURL r
  + \sum_{E\in\ET}\omega_{TE}\int_E (\bvec{v}_E\cdot \tangent_E) \, r
  \overset{\eqref{eq:def.SROT}}=
  -\int_T\SROT{T}\underline{\bvec{v}}_T \, r.\qedhere
  \]
\end{proof}

\subsection{Cohomology of the BGG complexes}
\label{sec:cohomology.bgg}

In this section, we derive the complexes obtained from the BGG diagram \eqref{eq:double.complex} and prove that their cohomologies are isomorphic to those of the corresponding continuous complexes.

The BGG complex derived from \eqref{eq:double.complex} is
\begin{equation}\label{derived-BGG}
  \begin{tikzcd}
    \text{DH($k+1$):}
    &[-2em] 0 \arrow{r}{}
    &[0em] \XSgrad{\Th} \arrow{r}{\Hess{h}}
    &[1.5em] \XdRrotS{\Th} \arrow{r}{\tROT{h}}
    &[1.5em] \Poly{k+1}(\Th)^2 \arrow{r}{}
    &[0em] 0
  \end{tikzcd}
\end{equation}
where
\begin{align*}
  \XdRrotS{\Th}\coloneq{}& \tXdRrot{\Th}\cap \Ker(\sskwh)\\
  ={}&
  \Big\{\underline{\bvec{\tau}}_h=((\bvec{\tau}_T)_{T\in\Th},(\bvec{\tau}_E)_{E\in\Eh})\,:\nonumber\\
  &\qquad\bvec{\tau}_T\in \Poly{k}(T)^{2\times 2}\text{ with $\bvec{\tau}(\bvec{x})\in\mathbb{S}$ for all $\bvec{x}\in T$}\quad\forall T\in\Th\,,\\
  &\qquad\bvec{\tau}_E\in\Poly{k+1}(E)^2\quad\forall E\in\Eh\Big\}
\end{align*}
and the discrete Hessian operator is given by $ \Hess{h}\coloneq\tGRAD{h}\circ \SGRAD{h}$.
The complex \eqref{derived-BGG} is referred to as the DH($k+1$) complex, standing for "Discrete Hessian" complex, where $k+1$ indicates the polynomial consistency degree of its differential operators, as proved in Lemma \ref{lem:pol.const.hess} for the discrete Hessian and \cite[(3.20)]{Di-Pietro.Droniou:23} for the discrete rotor (in the case of scalar complexes).
We refer the reader to \cite{Di-Pietro.Droniou.ea:26.1} for an application of the discrete complex \eqref{derived-BGG} to the design and analysis of a polytopal scheme for the Kirchhoff--Love plates problem.

\begin{lemma}[Commutation property and polynomial consistency of the discrete Hessian operator]
  \label{lem:pol.const.hess}
  For all $T\in\Th$,
  \begin{equation}\label{eq:com.prop.hess}
    \Hess{T}\ISgradT q =\IdRrotT \hess q \qquad \forall q \in \Cspace{2}(T),
  \end{equation}
  and the operator $\Hess{T}$ is consistent of degree $k+1$, that is,
  \begin{equation}\label{eq:pol.const.hess}
    \Prot{k+1}{T} \Hess{T}\ISgradT q = \hess q\qquad\forall q\in\Poly{k+3}(T),
  \end{equation}
  where $\Prot{k+1}{T}$ is the tensorised version of the potential on $\XdRrot{\Th}$ (built from the two-dimensional tangential trace of \cite[Eq.~(3.22), (3.23)]{Di-Pietro.Droniou:23} and transferred to the serendipity complex via \cite[Section 2]{Di-Pietro.Droniou:23*1}).
\end{lemma}
\begin{proof}
  Let $T\in\Th$. By the commutation property \eqref{eq:commutation.grad} of $\SGRAD{h}$ and the commutation property \cite[Eq.~(3.38)]{Di-Pietro.Droniou:23} of $\tGRAD{h}$, we have, for $q \in \Cspace{2}(\overline{T})$,
  \[
  \Hess{T}\ISgradT q = \tGRAD{T}\circ \SGRAD{T} \ISgradT q = \tGRAD{T} \IdRgradT (\GRAD q) = \IdRrotT (\GRAD\GRAD q)= \IdRrotT \hess q.
  \]
  Result \eqref{eq:pol.const.hess} follows by applying $\Prot{k+1}{T}$ to the above equality and using  its polynomial consistency (inferred from \cite[Proposition 3]{Di-Pietro.Droniou:23} and \cite[Proposition 7]{Di-Pietro.Droniou:23*1}).
\end{proof}

The twisted complex built  from \eqref{eq:double.complex} is
\begin{equation}\label{twisted-elasticity}
  \begin{tikzcd}[ampersand replacement=\&]
    0\arrow{r}\&[-1em]
    \begin{pmatrix}
      \XSgrad{\Th} \\
      \tXdRgrad{\Th}
    \end{pmatrix}
    \arrow{r}{
      \begin{pmatrix} \SGRAD{h} & -\Id \\ 0 &\tGRAD{h}  \end{pmatrix}
    }\&[3.5em]  \begin{pmatrix}
      \XSrot{\Th} \\
      \tXdRrot{\Th}
    \end{pmatrix} \arrow{r}{
      \begin{pmatrix}\SROT{h} &-\sskwh \\ 0 &\tROT{h}\end{pmatrix}
    } \&[3.8em] \begin{pmatrix}
      \Poly{k}(\Th)   \\
      \Poly{k+1}(\Th)^2
    \end{pmatrix} \arrow{r}{} \&[-1em] 0.
  \end{tikzcd}
\end{equation}
We assume that the domain is connected, possibly with holes, and
follow a dimension count argument to analyse the cohomology of \eqref{derived-BGG} and \eqref{twisted-elasticity}.

\begin{lemma}[Surjectivity of $\tROT{h}$]\label{lem:surjectivity.tROTh}
  In \eqref{derived-BGG},   $\tROT{h}: \XdRrotS{\Th}\to  \Poly{k+1}(\Th)^2$ is surjective.
\end{lemma}

\begin{proof}
  This result follows from a diagram chase on the diagrams in Figure \ref{fig:diag.chase} (a similar argument can be found in \cite{Arnold.Falk.ea:07}):
  For any $\bvec{w}_{h} \in \Poly{k+1}(\Th)^2$, we aim to find $\dofh{\bvec \tau} \in \XdRrotS{\Th}$ such that $\tROT{h}\dofh{\bvec \tau} = \bvec{w}_{h}$. To achieve this, we first find $\dofh{\bvec{\tilde \tau}} \in \tXdRrot{\Th}$ such that $\tROT{h} \dofh{\bvec{\tilde \tau}} = \bvec{w}_{h}$, which is possible because $\tROT{h}: \tXdRrot{\Th} \to \Poly{k+1}(\Th)^2$ is surjective \cite[Remark~10]{Di-Pietro.Droniou:23}. Then, we set $r_{h} \coloneqq  \sskwh \dofh{\bvec{\tilde{\tau}}} \in \Poly{k}(\Th)$. Using the surjectivity of $\SROT{h}: \XSrot{\Th} \to \Poly{k}(\Th)$, there exists $\dofh{\bvec{v}} \in \XSrot{\Th}$ such that $\SROT{h}\dofh{\bvec{v}} = r_{h}$. We then define $\dofh{\bvec \tau} \coloneqq \dofh{\bvec{\tilde \tau}}+ \tGRAD{h} \dofh{\bvec{v}} $. By the complex property, we have
  \[
  \tROT{h} \dofh{\bvec \tau} = \tROT{h} (\dofh{\bvec{\tilde \tau}}+ \tGRAD{h} \dofh{\bvec{v}}) = \tROT{h} \dofh{\bvec{\tilde \tau}} = \bvec{w}_{h},
  \]
  and, using the anti-commutativity of the diagram (Lemma \ref{lem:anticomm}), we find
  \[
  \sskwh \dofh{\bvec \tau} = \sskwh (\dofh{\bvec{\tilde \tau}}+ \tGRAD{h} \dofh{\bvec{v}}) = \sskwh \dofh{\bvec{\tilde \tau}} - \SROT{h}  \dofh{\bvec{v}} = r_{h} - r_{h} = 0.
  \]
  This implies that $\dofh{\bvec \tau} \in \XdRrotS{\Th}$, as required.
\end{proof}

\begin{figure}[h]
  \centering
  \begin{tikzcd}[column sep=2.5em]
    & \dofh{\bvec{\tilde \tau}} \arrow{r}{\tROT{h}} & \bvec{w}_{h}
    & \dofh{\bvec v} \arrow{r}{\SROT{h}} & r_{h} = \sskwh \dofh{\bvec{\tilde{\tau}}}
    \\
    & & \XSrot{\Th} \arrow{r}{\SROT{h}} & \Poly{k}(\Th) \\
    & \tXdRgrad{\Th} \arrow{r}{\tGRAD{h}} \arrow[leftrightarrow]{ur}{\Id} & \tXdRrot{\Th} \arrow{r}{\tROT{h}} \arrow{ur}{\sskwh} & \Poly{k+1}(\Th)^2 \\
    & & \dofh{\bvec \tau} = \dofh{\bvec{\tilde \tau}}+ \tGRAD{h} \dofh{\bvec{v}}  \arrow{r}{\tROT{h}} & \bvec{w}_{h}
  \end{tikzcd}
  \caption{Diagram chase in the proof of Lemma~\ref{lem:surjectivity.tROTh}: Step 1: Surjectivities of $\tROT{h}$ and $\SROT{h}$. Step 2: anti-commutativity. Step 3: Complex property.}
  \label{fig:diag.chase}
\end{figure}

\begin{lemma}[Kernel of the discrete Hessian operator]
  \label{lem:kernel.hess}
  The kernel of the discrete Hessian operator is spanned by the interpolates of affine functions on $\XSgrad{\Th}$, that is
  \begin{equation}
    \label{eq:kernel.hess}
    \Ker \Hess{h} = \ISgrad\Poly{1}(\Omega).
  \end{equation}
\end{lemma}

\begin{proof}
The inclusion $\supset$ follows from the commutation property \eqref{eq:com.prop.hess}. To prove the converse, let $\dofh{q}\in \XSgrad{\Th}$ be such that $\Hess{h}\dofh{q}=0$. Since $\Hess{h}\coloneq\tGRAD{h}\circ \SGRAD{h}$, we infer that  $\tGRAD{h} (\SGRAD{h}\dofh{q})=\dofh{\bvec{0}}$. Since DDR$(k+1)^2$ has the same cohomology as two copies of de Rham complex \cite[Theorem 1]{Di-Pietro.Droniou.ea:23}, and since $\Omega$ is connected, we have $\Ker\tGRAD{h}=\IdRgrad\Real^2$ and thus $\SGRAD{h}\dofh{q}=\IdRgrad \bvec{A}$ for some $\bvec{A}\in\Real^2$. We then deduce
\[
\SGRAD{h}\dofh{q} = \IdRgrad (\GRAD(\bvec{A}\cdot\bvec{x})) \overset{\eqref{eq:commutation.grad}}= \SGRAD{h}(\ISgrad (\bvec{A}\cdot\bvec{x})),
\]
so $\SGRAD{h}(\dofh{q}-\ISgrad (\bvec{A}\cdot\bvec{x}))=\dofh{\bvec{0}}$. Using then Theorem \ref{th:coho.stokes} and again the assumption that $\Omega$ is connected, we infer the existence of $c\in\Real$ such that $\dofh{q}-\ISgrad (\bvec{A}\cdot\bvec{x})=\ISgrad c$, which precisely shows that $\dofh{q}=\ISgrad (\bvec{A}\cdot\bvec{x}+c)$.
\end{proof}

Note that, by Lemma \ref{lem:kernel.hess}, $\dim \Ker (\Hess{h}) = 3 $. Moreover, one can easily use a dimension count to check that
\[
\dim \XSgrad{\Th} - \dim \XdRrotS{\Th} + \dim \Poly{k+1}(\Th)^2 = 3(\#\Vh-\#\Eh+\#\Th)= 3(-\beta_1 + 1),
\]
which is three times the Euler characteristic. Here, $\beta_1$ is the first Betti number (representing the number of holes), and we have used the relationship between the Euler characteristic $\chi=\#\Vh-\#\Eh+\#\Th$ and Betti numbers:
\[
\chi = \sum_{j=0}^{n} (-1)^{j} \beta_{j},
\]
along with the fact that $\beta_0 = 1$ since $\Omega$ is connected and $\beta_2 = 0$ since we are in dimension 2.
Note that, by Lemma \ref{lem:surjectivity.tROTh},
\[
\dim \Ker(\tROT{h}) = \dim \XdRrotS{\Th} - \dim \Poly{k+1}(\Th)^2,
\]
and
\[
\dim \Image(\Hess{h}) = \dim \XSgrad{\Th} - \dim \Ker (\Hess{h}) = \dim \XSgrad{\Th} -3.
\]
Therefore, the dimension of the cohomology group at the center of \eqref{derived-BGG} is
\[
\dim \Ker(\tROT{h}) - \dim \Image(\Hess{h}) = (\dim \XdRrotS{\Th} - \dim \Poly{k+1}(\Th)^2) - (\dim \XSgrad{\Th} - 3) = 3\beta_{1}.
\]
This implies that the dimensions of the cohomology groups of the discrete complex \eqref{derived-BGG} match those of the continuous one. Consequently, the discrete cohomology is isomorphic to its continuous counterpart.

The analysis of the cohomology of \eqref{twisted-elasticity} follows a similar dimension count argument.

\subsection{Comparison with finite element constructions}\label{sec:fe.comparison}

\begin{figure}[htbp]
  \centering
  \includegraphics[width=0.5\linewidth]{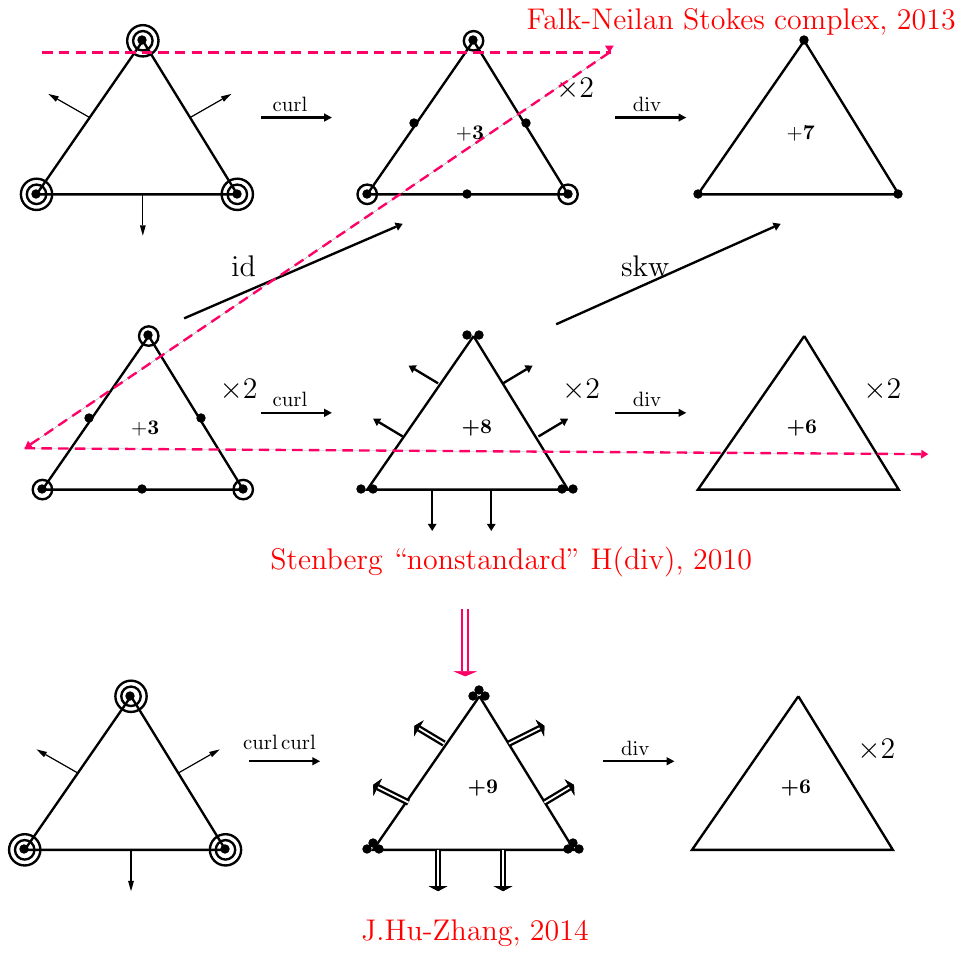}
  \caption{The BGG construction of the Hu--Zhang elasticity pair and the complex \cite{Christiansen.Hu.ea:18}. The first row is the Falk-Neilan Stokes complex. The second row is a de~Rham complex with enhanced vertex continuity. The $H(\mathrm{div})$ space is the ``nonstandard finite element'' by Stenberg \cite{Stenberg:10}.}
  \label{fig:hu-zhang}
\end{figure}

In this section, we compare the number of DOFs on triangles of the complexes built in the previous sections and corresponding finite element complexes from the BGG construction presented in \cite{Christiansen.Hu.ea:18}.
As illustrated in Figure \ref{fig:hu-zhang}, this construction is based on a diagram made of a Falk--Neilan Stokes complex (FN($\ell+1$)) \cite{Falk.Neilan:13} and a discrete de Rham complex \cite{Christiansen.Hu.ea:18}, resulting in the Hu--Zhang complex (HZ($\ell$)), which involves the Hu--Zhang element \cite{Hu.Zhang:15}. Here, as for DS($k$) and DDR$(k+1)^2$, the integers in the notations FN($\bullet$) and HZ($\bullet$) denote the degree of (minimal) polynomial exactness of all operators in the corresponding complex (we note that, for these finite element spaces, the degree of consistency decreases along the complex, so the one we consider here is that of the last differential operator).

\subsubsection{Description of the finite element spaces}

We review the definitions of the finite element spaces in Figure \ref{fig:hu-zhang} for completeness. We denote the diagram and the spaces therein by 
\begin{equation*}
  \begin{tikzcd}[column sep=3em, row sep=3em]
V^{0}_{h} \arrow{r}{\CURL} &V^{1}_{h}  \arrow{r}{\DIV}  & V^{2}_{h} \\
W^{0}_{h}  \arrow[ur, "\Id"]\arrow{r}{\CURL}& W^{1}_{h}\arrow{r}{\DIV} \arrow[ur, "\sskw"]& W^{2}_{h}.
  \end{tikzcd}
\end{equation*}

{\noindent --} $V^{0}_{h}$. The space $V^{0}_{h}$ is the Argyris finite element space. The local shape functions are quintic on a triangle. The degrees of freedom for $w\in V^0_T$ are
\begin{itemize}
\item $w(\bvec{x}_V), \partial_\ell w(\bvec{x}_V), \partial_\ell\partial_m w(\bvec{x}_V), \quad V\in \VT, ~\ell,m=1, 2~, m\le \ell$,
\item $\int_{E}\partial_{n}w\,ds, \quad E\in \ET$,
\end{itemize}
where $\partial_{n}w:=\normal\cdot \GRAD w$ is the normal derivative. An alternative definition of the global finite element space $V^{0}_{h}$ in terms of shape functions and inter-elements continuity is
$$
V^{0}_{h}:=\{u_{h}\in C^{1}(\Omega)\cap C^{2}(\Vh): u_{h}|_{T}\in \Poly{5}(T)\quad\forall T\in \mathcal{T}_{h}\},
$$
where $C^{1}(\Omega)\cap C^{2}(\Vh)$ denotes functions that are continuous on $\Omega$ and $C^{2}$ at all vertices of the triangulation $\Th$.
\medskip

{\noindent --} $V^{1}_{h}$. The space $V^{1}_{h}$ consists of two copies of the quartic Hermite finite element, one for each component of the vector field. The degrees of freedom for each component of $\bvec{w}=(w_1,w_2)\in V^1_T$ are
\begin{itemize}
\item $w_{i}(\bvec{x}_V), \partial_\ell w_{i}(\bvec{x}_V), \quad V\in \VT,~\ell=1,2$,
\item $\int_{E}w_{i}\,ds, \quad E\in \ET$,
\item $\int_{T}w_{i}q\,ds,  \quad q\in \Poly{1}(T)$.
\end{itemize}
An alternative definition of the global finite element space $V^{1}_{h}$ in terms of shape functions and interelement continuity is
$$
V^{1}_{h}:=\{\bvec{w}=(w_{1}, w_{2}): w_{i}\in C^{0}(\Omega)\cap C^{1}(\Vh): w_{i}|_{T}\in \Poly{4}(T)\quad\forall T\in \Th,\, \forall i=1, 2\}.
$$

{\noindent --} $V^{2}_{h}$. The space $V^{2}_{h}$ consists of piecewise cubic polynomials with vertex continuity. The degrees of freedom for $w\in V^2_T$ are
\begin{itemize}
\item $w(\bvec{x}_V), \quad V\in\VT$,
\item interior degrees of freedom.
\end{itemize}
An alternative definition of the global finite element space $V^{2}_{h}$ is
$$
V^{2}_{h}:=\{w\in  C^{0}(\mathcal V_h): w|_{T}\in \Poly{3}(T)\quad \forall T\in \Th\}.
$$

{\noindent --} $W^{0}_{h}$. To match the two rows precisely, we use the same finite element space at both ends of the identity map, i.e., we choose $W^{0}_{h}:=V^{1}_{h}$.\\
Since $W^{0}_{h}$ is a vector Lagrange finite element with enhanced vertex continuity, the remaining spaces in the bottom complex are the Whitney forms with the corresponding enhancement.
\medskip

{\noindent --} $W^{1}_{h}$. Two copies of the two-dimensional BDM finite element with vertex continuity, i.e., the Stenberg element  \cite{Stenberg:10}. Note that $W^{1}_{h}$ is a matrix-valued finite element. The shape function space consists of piecewise cubic polynomials for each component. The degrees of freedom of each component of $\bvec{w}=(w_{ij})_{i,j=1,2}\in W^1_T$ are
\begin{itemize}
\item $w_{ij}(\bvec{x}_V), \quad i, j=1, 2, ~V\in\VT$,
\item $\int_{E}(\bvec{w} \normal)_i q\,ds, \quad i=1, 2, ~ q\in \Poly{3}(E), ~ E\in \ET$, 
\item interior degrees of freedom.
\end{itemize}
\medskip

{\noindent --} $W^{2}_{h}$. Piecewise quadratic polynomials.

  With these choices of finite element spaces, the operator $\Id$ maps $W_h^0$ to $V_h^1$, and $\sskw$ maps $W_h^1$ to $V_h^2$. In the BGG reduction, one eliminates the spaces connected by the diagonal maps $\Id$ and $\sskw$. More precisely, $W_h^0$ and $V_h^1$ are completely removed since they coincide. Furthermore, eliminating $V_h^2$ from $W_h^1$ reduces the four vertex degrees of freedom of $W_h^1$ to three, and removes from its interior degrees of freedom those corresponding to $V_h^2$ (cf.~\cite{Christiansen.Hu.ea:18,Arnold.Falk.ea:06}; a rigorous argument involves defining a pseudo-inverse of $\sskw$, mapping from $V_h^2$ to $W_h^1$, by setting all remaining degrees of freedom in $W_h^1$ to zero). One then connects the two rows in the diagram via a zig-zag construction to obtain the Hu--Zhang elasticity complex:
\begin{equation*}
  \begin{tikzcd}[column sep=3em]
V^{0}_{h} \arrow[r, "{\CURL}"] & \Sigma_{h} \arrow[r, "{\DIV}"] & W^{2}_{h}.
  \end{tikzcd}
\end{equation*}
Here the stress space is defined by
\[
\Sigma_{h}:=\{\sigma_h:\Omega\to\mathbb{S}\,:\,(\sigma_h)_{|T}\in \Poly{3}(T)^{2\times 2} \quad\forall T\in\Th\,,\;
 \sigma_h \cdot \bvec{n} \text{ is single-valued}, \;\sigma_h\in C^0(\mathcal{V}_h)\}.
\]
In terms of shape functions and degrees of freedom, the space $\Sigma_{h}$ consists of symmetric matrix-valued functions whose components are piecewise cubic polynomials. The degrees of freedom are given by
\begin{itemize}
\item $w_{ij}(\bvec{x}_V), \quad i,j=1,2,\ i \le j,\ V \in \VT$,
\item $\int_{E} \bvec{w}_{i} \cdot \normal\, q \, ds, \quad i=1,2,\ q \in \Poly{1}(E),\ E\in\ET$,
\item interior degrees of freedom.
\end{itemize}

\subsubsection{Comparison between finite element and polytopal complexes}

To ensure a meaningful comparison, we consider complexes with the same degree of polynomial accuracy. The comparison between DS($k$) and FN($k$) is given in Table \ref{tab:comparison.dofs.stokes}, while Table \ref{tab:comparison.dofs.hess} concerns DH($k+1$) and HZ($k+1$).

\begin{table}[h]
  \centering
  \adjustbox{width=\textwidth}{
    \begin{tabular}{c|c|c|c|c}
      \toprule
      \multicolumn{2}{c|}{} & \multicolumn{3}{c}{Continuous space}\\
      \midrule
      Discrete complex &   & $H_2(\Omega)$ & $H_1(\Omega)^2$ & $L^2(\Omega)$ \\
      \midrule
      \multirow{4}{*}{DS$(k)$ ($k\geq 0$)}
      & \textbf{Total DOFs per triangle} & $12 + \frac{1}{2}(11k+k^2)$  & $12 +7k+k^2$ &  $\frac{1}{2}(k+2)(k+1)$\\
      \cmidrule(lr){2-5}
      & per vertex   & 3 & 2 & 0\\
      & per edge & $2k+1$ & $2(k+1)$ & 0 \\
      & in the element & $\frac{1}{2}(k-1)k$ & $k(k+1)$ & $\frac{1}{2}(k+2)(k+1) $\\
      \midrule
      \multirow{4}{*}{FN($k$) ($k\geq 3$)}
      & \textbf{Total DOFs per triangle} & $6 + \frac{1}{2}(7k+k^2)$  & $6 + 5k+k^2$ & $\frac{1}{2}(k+2)(k+1)$ \\
      \cmidrule(lr){2-5}
      & per vertex  & 6 & 6 &  1\\
      & per edge  &  $2k-5$ & $2(k-2)$ & 0\\
      & in the element & $\frac{1}{2}(k-3)(k-2)$  & $(k-1)k$ & $\frac{1}{2}(k+2)(k+1)-3$ \\
      \bottomrule
    \end{tabular}
  }
  \caption{Comparison of DOFs in the discrete Stokes complex, for each space, on a triangle.}
  \label{tab:comparison.dofs.stokes}
\end{table}

\begin{table}[h]
  \centering
  \adjustbox{width=\textwidth}{
    \begin{tabular}{c|c|c|c|c}
      \toprule
      \multicolumn{2}{c|}{} & \multicolumn{3}{c}{Continuous space}\\
      \midrule
      Discrete complex &  & $H_2(\Omega)$ & $\boldsymbol H_{\boldsymbol \ROT}(\Omega, \mathbb{S})$ & $L^2(\Omega)^2$ \\
      \midrule
      \multirow{4}{*}{DH$(k+1)$ ($k\geq 0$)}
      & \textbf{Total DOFs per triangle} & $12 + \frac{1}{2}(11k+k^2)$
      & $15 +\frac{3}{2}(7k+k^2)$ &   $(k+2)(k+3)$ \\
      \cmidrule(lr){2-5}
      & per vertex   & 3 & 0 & 0\\
      & per edge & $2k+1$ & $2k+4$ & 0 \\
      & in the element & $\frac{1}{2}(k-1)k$ & $\frac{3}{2}(k+1)(k+2)$ &  $(k+2)(k+3)$ \\
      \midrule
      \multirow{4}{*}{HZ($k+1$) ($k\geq 1$)}
      & \textbf{Total DOFs per triangle} & $15 + \frac{1}{2}(11k+k^2)$
      & $18 +\frac{3}{2}(7k+k^2)$
      & $(k+2)(k+3)$ \\
      \cmidrule(lr){2-5}
      & per vertex  & 6 & 3 &  0\\
      & per edge  &  $2k-1$ & $2k+2$ & 0\\
      & in the element & $\frac{1}{2}(k-1)k$  & $\frac{3}{2}(k+1)(k+2)$ & $(k+2)(k+3)$ \\
      \bottomrule
    \end{tabular}
  }
  \caption{Comparison of DOFs in the discrete Hessian complex, for each space, on a triangle.}
  \label{tab:comparison.dofs.hess}
\end{table}

These tables show that DS($k$) is slightly more expensive than the corresponding Falk--Neilan complex.
The spaces in DH($k+1$) have, on the other hand, fewer degrees of freedom than their counterparts in HZ($k+1$).
We notice that the dimension of the spaces in the polytopal complexes could be further reduced by taking full advantage of serendipity in the spirit of \cite{Di-Pietro.Droniou:23*1,Botti.Di-Pietro.ea:23}.
This topic will be addressed in a forthcoming work.

We also notice that, besides being applicable on generic polygonal meshes (which can lead to more efficient meshing of complicated domains than by triangles), the complexes we designed also provide lower-order (and thus cheaper) versions than those accessible via finite element constructions.
For example, the spaces of DS($0$) have 12/12/1 degrees of freedom per triangle while the spaces in the lowest-order Falk--Neilan complex have 21/30/10 degrees of freedom. For the Hessian complexes, DH($1$) has 12/15/6 DOFs per triangle while HZ($2$) has 21/18/12 DOFs.
In some circumstances, a lower-order method can be preferred over a higher-order method, for example in the case of non-linear problems, and/or when the solution is not expected to be smooth (and, thus, the more expensive higher-order method cannot realise an improved accuracy).
When solving systems using $p$-multigrid, being able to go lower in the degree of the method can also be beneficial.

Figure \ref{fig:double.comple.k.zero} provides a representation of the DOFs in the lowest-order BGG diagram \eqref{eq:double.complex}, equivalent on a hexagonal element to the lowest order diagram of \cite{Christiansen.Hu.ea:18} represented in Figure \ref{fig:hu-zhang}.

\begin{figure}[htbp]
  \centering
  \begin{equation*}
    \begin{tikzcd}[column sep=1.5em]
      \text{DS(0):}
      & 0 \arrow{r}{}
      & \DrawHtwo
      \arrow{r}{\underline{\bvec{G}}^{0}_{2,T}}
      &[2em] \DrawHone \times 2 \arrow{r}{R^0_{1,T}}
      &[2em]\DrawPzero\arrow{r}{}
      & 0 \\
      \text{DDR$(1)^2$:}
      & 0 \arrow{r}{}
      & \DrawHone \times 2 \arrow{r}{\underline{\bvec G}_{1,T}^{1}}\arrow[leftrightarrow]{ur}{\Id}
      & \DrawHrot \times 2 \arrow{r}{\bvec{R}_{\VROT,T}^{1}}\arrow{ur}{\sskw\nolimits_T}
      & \DrawPone \times 2 \arrow{r}{}
      & 0.
    \end{tikzcd}
  \end{equation*}
  \caption{Schematic representations of DOFs for the lowest-order BGG diagram \eqref{eq:double.complex}.\label{fig:double.comple.k.zero}
  }
\end{figure}

\section*{Acknowledgements}

The authors acknowledge the funding of the European Union via the ERC Synergy, NEMESIS, project number 101115663.

KH was supported in part by the ERC Starting project GeoFEM, project number 101164551, and a Royal Society University Research Fellowship (URF$\backslash$R1$\backslash$221398). Part of the work was carried out during KH's visit at University of Montpellier supported by the NEMESIS project. KH is grateful for the hospitality and the support from NEMESIS.

Views and opinions expressed are however those of the authors only and do not necessarily reflect those of the European Union or the European Research Council Executive Agency. Neither the European Union nor the granting authority can be held responsible for them.

\appendix

\section{Notation}
\begin{center}
  \renewcommand{\arraystretch}{1.2}
  \begin{longtable}{ccc}
    \toprule
    \textbf{Symbol} & \textbf{Continuous counterpart} & \textbf{Definition} \\
    \midrule
    \multicolumn{3}{c}{Discrete spaces} \\
    \midrule
    $\XSgrad{\Th}$
    & $H_2(\Omega)$
    & \eqref{eq:XSgrad} \\
    $\XSrot{\Th}$
    & $H_1(\Omega)^2$
    & \eqref{eq:def.Xgrad} \\
    $\tXdRrot{\Th}$
    & $\boldsymbol{H}_{\VROT}(\Omega)$
    & \eqref{eq:def.Xrot} \\
    \midrule
    \multicolumn{3}{c}{Interpolators} \\
    \midrule
    $\ISgrad$
    & ---
    & \eqref{eq:def.ISgrad} \\
    $\IdRgrad$
    & ---
    & \eqref{eq:def.IdRgrad} \\
    $\IdRrot$
    & ---
    & \eqref{eq:def.IdRrot} \\
    \midrule
    \multicolumn{3}{c}{Discrete differential operators} \\
    \midrule
    $\SGRAD{h}$
    & $\GRAD : H_2(\Omega) \to H_1(\Omega)$
    & \eqref{eq:def.nablah} \\
    $\SROT{h}$
    & $\ROT$
    & \eqref{eq:SROT.h} \\
    $\tGRAD{h}$
    & $\GRAD : H_1(\Omega)^2 \to \boldsymbol{H}_{\VROT}(\Omega)$
    & \eqref{eq:def.uGh1} \\
    $\tROT{h}$
    & $\VROT$
    & \eqref{eq:tROT.h} \\
    \midrule
    \multicolumn{3}{c}{Algebraic operators} \\
    \midrule
    $\sskw$
    & ---
    & \eqref{eq:sskw} \\
    $\sskwh$
    & $\sskw$
    & \eqref{eq:def.sskwh} \\
    \bottomrule
  \end{longtable}
\end{center}


\printbibliography

@InProceedings{   Arnold.Falk.ea:06,
  title         = {{Differential complexes and stability of finite element
                  methods II: The elasticity complex}},
  author        = {Arnold, Douglas N and Falk, Richard S and Winther,
                  Ragnar},
  booktitle     = {Compatible spatial discretizations},
  pages         = {47--67},
  year          = {2006},
  organization  = {Springer}
}

@Article{         Arnold.Falk.ea:07,
  title         = {Mixed finite element methods for linear elasticity with
                  weakly imposed symmetry},
  author        = {Arnold, Douglas and Falk, Richard and Winther, Ragnar},
  journal       = {Math. Comp.},
  volume        = {76},
  number        = {260},
  pages         = {1699--1723},
  year          = {2007}
}

@Article{         Arnold.Falk.ea:10,
  author        = {Arnold, Douglas N. and Falk, Richard S. and Winther,
                  Ragnar},
  title         = {Finite element exterior calculus: from {H}odge theory to
                  numerical stability},
  journal       = {Bull. Amer. Math. Soc. (N.S.)},
  volume        = {47},
  year          = {2010},
  number        = {2},
  pages         = {281--354},
  doi           = {10.1090/S0273-0979-10-01278-4}
}

@Article{         Arnold.Hu:21,
  title         = {Complexes from complexes},
  author        = {Arnold, Douglas N and Hu, Kaibo},
  journal       = {Found. Comput. Math.},
  pages         = {1--36},
  year          = {2021},
  publisher     = {Springer}
}

@Article{         Arnold.Qin:92,
  title         = {{Quadratic velocity/linear pressure Stokes elements}},
  author        = {Arnold, Douglas N and Qin, Jinshui},
  journal       = {Advances in computer methods for partial differential
                  equations},
  volume        = {7},
  pages         = {28--34},
  year          = {1992}
}

@Article{         Arnold.Winther:02,
  title         = {Mixed finite elements for elasticity},
  author        = {Arnold, Douglas N and Winther, Ragnar},
  journal       = {Numer. Math.},
  volume        = {92},
  pages         = {401--419},
  year          = {2002},
  publisher     = {Springer}
}

@Book{            Arnold:18,
  author        = {Arnold, D.},
  title         = {Finite Element Exterior Calculus},
  publisher     = {SIAM},
  year          = {2018},
  doi           = {10.1137/1.9781611975543}
}

@Article{         Beirao-da-Veiga.Mora.ea:19,
  author        = {Beir{\~a}o da Veiga, L. and Mora, D. and Vacca, G.},
  title         = {The {S}tokes complex for virtual elements with application
                  to {N}avier-{S}tokes flows},
  journal       = {J. Sci. Comput.},
  volume        = {81},
  year          = {2019},
  number        = {2},
  pages         = {990--1018},
  doi           = {10.1007/s10915-019-01049-3}
}

@Article{         Bernstein.Gelfand.ea:75,
  title         = {{Differential operators on the base affine space and a
                  study of g-modules}},
  author        = {Bernstein, IN and Gelfand, Israel M and Gelfand, Sergei
                  I},
  journal       = {Lie groups and their representations (Proc. Summer School,
                  Bolyai J{\'a}nos Math. Soc., Budapest, 1971)},
  pages         = {21--64},
  year          = {1975}
}

@Article{         Bonaldi.Di-Pietro.ea:25,
  title         = {An exterior calculus framework for polytopal methods},
  author        = {Bonaldi, F. and Di Pietro, D. A. and Droniou, J. and Hu,
                  K.},
  journal       = {J. Eur. Math. Soc.},
  year          = {2025},
  doi           = {10.4171/JEMS/1602},
  note          = {Published online}
}

@Article{         Bonizzoni.Hu.ea:25,
  author        = {Bonizzoni, Francesca and Hu, Kaibo and Kanschat, Guido and
                  Sap, Duygu},
  title         = {Discrete tensor product {BGG} sequences: splines and
                  finite elements},
  journal       = {Math. Comp.},
  volume        = {94},
  year          = {2025},
  number        = {352},
  pages         = {517--549},
  doi           = {10.1090/mcom/3969}
}

@Article{         Boon.Duran.ea:25,
  title         = {{Mixed finite element methods for linear Cosserat
                  equations}},
  author        = {Boon, Wietse Marijn and Duran, Omar and Nordbotten, Jan
                  Martin},
  journal       = {SIAM J. Numer. Anal.},
  volume        = {63},
  number        = {1},
  pages         = {306--333},
  year          = {2025},
  publisher     = {SIAM}
}

@Article{         Botti.Di-Pietro.ea:23,
  author        = {Botti, Michele and Di Pietro, Daniele A. and Salah,
                  Marwa},
  title         = {A serendipity fully discrete div-div complex on polygonal
                  meshes},
  journal       = {Comptes Rendus Mécanique},
  year          = { 2023},
  volume        = {351(S1)},
  doi           = {10.5802/crmeca.150}
}

@Article{         Cap.Hu:24,
  title         = {{BGG sequences with weak regularity and applications}},
  author        = {{\v{C}}ap, Andreas and Hu, Kaibo},
  journal       = {Found. Comput. Math.},
  volume        = {24},
  number        = {4},
  pages         = {1145--1184},
  year          = {2024},
  publisher     = {Springer}
}

@Article{         Cap.Slovak.ea:01,
  title         = {{Bernstein-Gelfand-Gelfand sequences}},
  author        = {{\v{C}}ap, Andreas and Slov{\'a}k, Jan and Sou{\v{c}}ek,
                  Vladim{\'\i}r},
  journal       = {Annals of Mathematics},
  volume        = {154},
  number        = {1},
  pages         = {97--113},
  year          = {2001},
  publisher     = {JSTOR}
}

@InCollection{    Chen.Huang:22,
  title         = {{Discrete Hessian complexes in three dimensions}},
  author        = {Chen, Long and Huang, Xuehai},
  booktitle     = {The virtual element method and its applications},
  pages         = {93--135},
  year          = {2022},
  publisher     = {Springer}
}

@Article{         Chen.Huang:22*1,
  title         = {Finite elements for div-and divdiv-conforming symmetric
                  tensors in arbitrary dimension},
  author        = {Chen, Long and Huang, Xuehai},
  journal       = {SIAM Journal on Numerical Analysis},
  volume        = {60},
  number        = {4},
  pages         = {1932--1961},
  year          = {2022},
  publisher     = {SIAM}
}

@Article{         Chen.Huang:22*2,
  title         = {A finite element elasticity complex in three dimensions},
  author        = {Chen, Long and Huang, Xuehai},
  journal       = {Mathematics of Computation},
  volume        = {91},
  number        = {337},
  pages         = {2095--2127},
  year          = {2022}
}

@Article{         Chen.Huang:24,
  title         = {Finite element complexes in two dimensions},
  author        = {Chen, Long and Huang, Xuehai},
  journal       = {Scientia Sinica Mathematica},
  volume        = {55},
  number        = {8},
  pages         = {1593--1626},
  year          = {2025},
  doi           = {10.1360/SCM-2023-0169}
}

@Article{         Christiansen.Hu.ea:18,
  title         = {{Nodal finite element de Rham complexes}},
  author        = {Christiansen, Snorre H and Hu, Jun and Hu, Kaibo},
  journal       = {Numer. Math.},
  volume        = {139},
  number        = {2},
  pages         = {411--446},
  year          = {2018},
  publisher     = {Springer}
}

@Article{         Christiansen.Hu:18,
  title         = {{Generalized finite element systems for smooth
                  differential forms and Stokes’ problem}},
  author        = {Christiansen, Snorre H and Hu, Kaibo},
  journal       = {Numer. Math.},
  volume        = {140},
  pages         = {327--371},
  year          = {2018},
  publisher     = {Springer}
}

@Article{         Christiansen.Hu:23,
  title         = {Finite element systems for vector bundles: elasticity and
                  curvature},
  author        = {Christiansen, Snorre H and Hu, Kaibo},
  journal       = {Found. Comput. Math.},
  volume        = {23},
  number        = {2},
  pages         = {545--596},
  year          = {2023},
  publisher     = {Springer}
}

@Article{         Di-Pietro.Droniou.ea:20,
  author        = {Di Pietro, D. A. and Droniou, J. and Rapetti, F.},
  title         = {Fully discrete polynomial {de Rham} sequences of arbitrary
                  degree on polygons and polyhedra},
  journal       = {Math. Models Methods Appl. Sci.},
  year          = {2020},
  volume        = {30},
  number        = {9},
  pages         = {1809-1855},
  doi           = {10.1142/S0218202520500372}
}

@Article{         Di-Pietro.Droniou.ea:23,
  title         = {Cohomology of the discrete de Rham complex on domains of
                  general topology},
  volume        = {60},
  doi           = {10.1007/s10092-023-00523-7},
  number        = {2},
  journal       = {Calcolo},
  publisher     = {Springer Science and Business Media LLC},
  author        = {Di Pietro, Daniele A. and Droniou, Jérôme and Pitassi,
                  Silvano},
  year          = {2023},
  month         = may
}

@Article{         Di-Pietro.Droniou.ea:26.1,
  author        = {Di Pietro, D. A. and Droniou, J. and Hu, K. and Leroy,
                  A.},
  title         = {Analytical properties of polygonal Stokes and BGG Hessian
complexes with application to Kirchhoff–Love plates},
  year          = {2026},
  note          = {Submitted}
}

@Book{            Di-Pietro.Droniou:20,
  author        = {Di Pietro, D. A. and Droniou, J.},
  title         = {The {Hybrid High-Order} method for polytopal meshes},
  subtitle      = {Design, analysis, and applications},
  publisher     = {Springer International Publishing},
  year          = {2020},
  series        = {Modeling, Simulation and Application},
  number        = {19},
  doi           = {10.1007/978-3-030-37203-3}
}

@Article{         Di-Pietro.Droniou:21,
  author        = {Di Pietro, D. A. and Droniou, J.},
  title         = {An arbitrary-order method for magnetostatics on polyhedral
                  meshes based on a discrete {de Rham} sequence},
  year          = {2021},
  journal       = {J. Comput. Phys.},
  volume        = {429},
  number        = {109991},
  doi           = {10.1016/j.jcp.2020.109991}
}

@Article{         Di-Pietro.Droniou:23,
  author        = {Di Pietro, Daniele A. and Droniou, J\'er\^ome},
  title         = {An arbitrary-order discrete de Rham complex on polyhedral
                  meshes: Exactness, Poincar\'e inequalities, and
                  consistency},
  journal       = {Found. Comput. Math.},
  volume        = {23},
  pages         = {85--164},
  year          = {2023},
  doi           = {10.1007/s10208-021-09542-8}
}

@Article{         Di-Pietro.Droniou:23*1,
  author        = {Di Pietro, Daniele A. and Droniou, J\'er\^ome},
  title         = {Homological- and analytical-preserving serendipity
                  framework for polytopal complexes, with application to the
                  DDR method},
  pages         = {191--225},
  year          = {2023},
  volume        = {57},
  issue         = {1},
  journal       = {{M2AN} Math. Model. Numer. Anal.},
  doi           = {10.1051/m2an/2022067}
}

@Article{         Di-Pietro:24,
  title         = {An arbitrary-order discrete rot-rot complex on polygonal
                  meshes with application to a quad-rot problem},
  author        = {Di Pietro, D. A.},
  year          = {2024},
  journal       = {IMA J. Numer. Anal.},
  volume        = {44},
  number        = {3},
  pages         = {1699--1730},
  doi           = {10.1093/imanum/drad045}
}

@Article{         Falk.Neilan:13,
  title         = {{Stokes complexes and the construction of stable finite
                  elements with pointwise mass conservation}},
  author        = {Falk, Richard S and Neilan, Michael},
  journal       = {SIAM J. Numer. Anal.},
  volume        = {51},
  number        = {2},
  pages         = {1308--1326},
  year          = {2013},
  publisher     = {SIAM}
}

@Article{         Hanot:23,
  author        = {Hanot, Marien-Lorenzo},
  title         = {An arbitrary-order fully discrete {S}tokes complex on
                  general polyhedral meshes},
  journal       = {Math. Comp.},
  volume        = {92},
  year          = {2023},
  number        = {343},
  pages         = {1977--2023},
  doi           = {10.1090/mcom/3837}
}

@Article{         Hu.Liang.ea:22,
  title         = {{Conforming Finite Element DIVDIV Complexes and the
                  Application for the Linearized Einstein--Bianchi System}},
  author        = {Hu, Jun and Liang, Yizhou and Ma, Rui},
  journal       = {SIAM Journal on Numerical Analysis},
  volume        = {60},
  number        = {3},
  pages         = {1307--1330},
  year          = {2022},
  publisher     = {SIAM}
}

@Article{         Hu.Lin.ea:23,
  title         = {Finite elements for symmetric and traceless tensors in
                  three dimensions},
  author        = {Hu, Kaibo and Lin, Ting and Shi, Bowen},
  journal       = {arXiv preprint arXiv:2311.16077},
  year          = {2023}
}

@Article{         Hu.Lin.ea:25,
  author        = {Hu, Kaibo and Lin, Ting and Zhang, Qian},
  title         = {Distributional {H}essian and divdiv complexes on
                  triangulation and cohomology},
  journal       = {SIAM J. Appl. Algebra Geom.},
  volume        = {9},
  year          = {2025},
  number        = {1},
  pages         = {108--153},
  doi           = {10.1137/23M1623860}
}

@Article{         Hu.Zhang:15,
  title         = {A family of conforming mixed finite elements for linear
                  elasticity on triangular grids},
  author        = {Hu, Jun and Zhang, Shangyou},
  journal       = {Science China Mathematics},
  volume        = {58},
  pages         = {297--307},
  doi           = {10.1007/s11425-014-4953-5},
  year          = {2015}
}

@InCollection{    Hu:22,
  title         = {Nonlinear elasticity complex and a finite element diagram
                  chase},
  author        = {Hu, Kaibo},
  booktitle     = {Approximation Theory and Numerical Analysis Meet Algebra,
                  Geometry, Topology},
  series        = {Springer INdAM Series},
  volume        = {60},
  pages         = {231--252},
  year          = {2024},
  publisher     = {Springer},
  doi           = {10.1007/978-981-97-6508-9_11}
}

@Article{         Johnson.Mercier:78,
  title         = {Some equilibrium finite element methods for
                  two-dimensional elasticity problems},
  author        = {Johnson, Claes and Mercier, Bertrand},
  journal       = {Numer. Math.},
  volume        = {30},
  pages         = {103--116},
  year          = {1978},
  publisher     = {Springer}
}

@Book{            Lai.Schumaker:07,
  title         = {Spline functions on triangulations},
  author        = {Lai, Ming-Jun and Schumaker, Larry L},
  number        = {110},
  year          = {2007},
  publisher     = {Cambridge University Press}
}

@Article{         Stenberg:10,
  title         = {A nonstandard mixed finite element family},
  author        = {Stenberg, Rolf},
  journal       = {Numer. Math.},
  volume        = {115},
  pages         = {131--139},
  year          = {2010},
  publisher     = {Springer}
}


\end{document}